\documentclass[default]{sn-jnl}

\jyear{2024}%
\graphicspath{{Figures/}}
\usepackage{bm}
\usepackage{xcolor}
\usepackage{hyperref}
\usepackage{subfigure}
\usepackage{colortbl}
\usepackage{graphicx}
\usepackage{amsmath}
\usepackage{amssymb}
\usepackage{amsthm}
\usepackage{listings}
\usepackage{xcolor}
\usepackage{multirow}

\newtheorem{remark}{Remark}[section]

\begin{document}

\title[Quasi-interpolation projectors for Subdivision Surfaces]{Quasi-interpolation projectors for Subdivision Surfaces}

\author{\fnm{Hailun} \sur{Xu}}
\author*{\fnm{Hongmei} \sur{Kang}}
\email{khm@suda.edu.cn}


\affil{\orgdiv{School of Mathematical Sciences}, \orgname{Soochow University}, \orgaddress{\street{No.1 Road Shizi}, \city{Suzhou}, \postcode{215006}, \state{Jiangsu}, \country{China}}}




\abstract{
Subdivision surfaces are considered as an extension of splines to accommodate models with complex topologies, making them useful for addressing PDEs on models with complex topologies in isogeometric analysis. This has generated a lot of interest in the field of subdivision space approximation. The quasi-interpolation offers a highly efficient approach for spline approximation, eliminating the necessity of solving large linear systems of equations. Nevertheless, the lack of analytical expressions at extraordinary points on subdivision surfaces makes traditional techniques for creating B-spline quasi-interpolants inappropriate for subdivision spaces. To address this obstacle, this paper innovatively reframes the evaluation issue associated with subdivision surfaces as a correlation between subdivision matrices and limit points, offering a thorough method for quasi-interpolation specifically designed for subdivision surfaces. This developed quasi-interpolant, termed the subdivision space projection operator, accurately reproduces the subdivision space. We provide explicit quasi-interpolation formulas for various typical subdivision schemes. Numerical experiments demonstrate that the quasi-interpolants for Catmull-Clark and Loop subdivision exhibit third-order approximation in the \(L_2\) norm and second-order in the \(L_\infty\) norm. Furthermore, the modified Loop subdivision quasi-interpolant achieves optimal approximation rates in both the \(L_2\) and \(L_\infty\) norms.

}




\maketitle
\section{Introduction}

The subdivision surface technique provides a simple and effective algorithm for modeling and processing arbitrary topology smooth surfaces, which have become powerful tools for geometric modeling and computer graphics \cite{Derose,dyn_levin_2002,liao2017subdivision,peters2008subdivision}. Subdivision defines a smooth curve or surface as the limit of a sequence of successive refinements. Each iteration generates a smoother and more refined mesh. The limit surface is globally smooth. The well-known Doo-Sabin \cite{DOO1978356} and Catmull-Clark subdivisions \cite{CATMULL1978350}  degenerate into quadratic and cubic B-splines when there are no extraordinary points, respectively. The triangle-based Loop subdivision \cite{Loop} is a generalization of the quartic $C^2$-continuous box splines. Subdivision surfaces are currently utilized in isogeometric analysis to handle partial differential equations (PDEs) on models with complex topologies \cite{Thinshell,CIRAK2002137,KANGmodifiedloop,igasubdthesis,XIE2020101867}, demonstrating promising prospects in the finite element analysis of thin shells.

Approximation is a fundamental yet crucial numerical problem, especially within the context of subdivision surfaces. In \cite{Halstead1993},  an efficient interpolation method with Catmull-Clark subdivision is proposed. The resulting subdivision surface interpolates the vertices of any topological type of mesh while maintaining fairness. For the data obtained through CT scanning, the Catmull-Clark subdivision surfaces are fitted by the least-squares method in \cite{WeiyinMa}. For dense triangle meshes, a topology- and feature-preserving mesh simplification algorithm is employed to obtain a coarse mesh firstly, and then the control mesh of Loop subdivision surface is fitted by the least squares method in \cite{Ma2002,MA2004525}. However, these methods all need to solve large linear systems and are very time-consuming.

Quasi-interpolation is a general method for constructing approximations at a low computational cost, which was originally proposed to avoid solving linear equations globally, and has good numerical stability and approximation ability.
Studies on the quasi-interpolation of B-splines date back to the last century \cite{localapprox,quasimulti,deBoor73}. In \cite{examples}, a more practical recipe for deriving local spline projectors using quasi-interpolation for B-spline spaces is proposed. The concise quasi-interpolation formula for B-splines greatly simplifies the construction of local spline projectors in hierarchical spline spaces  \cite{hierarquasi,effortless,thbproj,LRQI} and analysis-suitable T-splines \cite{KANG2022102147}. For splines on triangulation, box spline projectors are discussed in \cite{LYCHE2008416} based on local interpolation and local $L_2$-norm inner product projection, respectively.

Although the subdivision surfaces are viewed as an extension of B-splines and box splines to models with complex topologies, the presence of extraordinary points poses a challenge in creating quasi-interpolation of subdivision surfaces because of the absence of analytical expressions at these points. In this paper, we still follow the general recipe of constructing quasi-interpolation for various splines \cite{examples} but with a new technique for addressing the local interpolation challenge. 
We construct the local interpolation problem based on the subdivision matrix and limit position, overcoming the cumbersome evaluation of the subdivision basis functions and the difficulty associated with deriving explicit solutions to the local interpolation problem. This proposed quasi-interpolation is very simple and intuitive, and easy to construct quasi-interpolation for general subdivision schemes.

In \cite{2001Approximation}, a quasi-interpolation method for Loop subdivision was developed to assist in proving approximation estimates. This quasi-interpolation operator involves complex integration and averaging calculations. Another study \cite{2001Fitting} presented a heuristic approach to construct the quasi-interpolation for the Catmull-Clark subdivision. However, this method can not reproduce the subdivision spaces and offers a lower approximation order.
The proposed quasi-interpolation in this paper is a subdivision space projector, thus it can reproduce the subdivision space.
Numerical experiments demonstrate that the quasi-interpolants of Loop and Catmull–Clark subdivisions achieve an approximation order of $3$ in the $L_2$ norm and $2$ in the $L_\infty$ norm, consistent with the approximation estimates presented in \cite{2001Approximation}. Furthermore, we investigate the quasi-interpolation of the modified Loop subdivision method introduced by \cite{KANGmodifiedloop}, which is designed to tackle the suboptimal convergence problem in Loop subdivision.
Through the reduction of subdominant eigenvalues at extraordinary points, the quasi-interpolation of the modified Loop subdivision method can attain the same optimal approximation order as in regular domains.

This paper is organized as follows. In Section \ref{s.generalrecipe}, we present the general recipe of constructing subdivision surface projectors based on the subdivision matrix and limit position matrix. In Section \ref{examples}, we provide the explicit expressions of the quasi-interpolation of Catmull-Clark, Loop, and modified Loop subdivisions, respectively. In Section \ref{Sec4}, we give several strategies for processing boundary vertices and vertices that are adjacent to extraordinary points. In Section \ref{nums}, numerical experiments are conducted to demonstrate the approximation orders of the proposed quasi-interpolations. Finally, Section \ref{conclu} summarizes the paper and discusses future work.

\section{Quasi-interpolation for Subdivision surfaces \label{s.generalrecipe}}
For any given control net, we can induce a subdivision surface $\Omega$ and the subdivision space  $\mathcal{S}$ defined on it, which contains a set of linearly independent subdivision basis functions, $\mathcal{S}=span\{B_1,B_2,\cdots,B_n\}$.
For a function $f \in C\left(\Omega\right)$, the quasi-interpolation $Q:C\left(\Omega\right)\rightarrow \mathcal{S}$ in the subdivision space is defined by
\begin{equation}
	\label{QIsB}
	Q\left(f\right)=\sum_{i=1}^n\lambda_{i}(f)B_i(x),~~x\in\Omega,
\end{equation}
where $\lambda_{i}$ is a linear functional and typically defined in the form of point functional,
\begin{equation}
	\label{lamdaK}
	\lambda_{i}(f)=\sum_k\omega_{i,k}f(x_{i,k}),~~\omega_{i,k}\in\mathbb{R}.
\end{equation}
To obtain a subdivision space projector, the quasi-interpolation $Q(f)$ is required to reproduce the whole subdivision space $\mathcal{S}$, that is 
\[Q(f)=f,~~\forall f \in \mathcal{S}.\]

The general procedure of determining quasi-interpolation projectors for B-spline spaces proposed in \cite{examples} has been widely used for locally refinable splines, such as hierarchical splines \cite{effortless,hierarquasi}, THB-splines\cite{thbproj}, AS T-splines \cite{dualast}, etc. We adopt the same procedure for constructing subdivision space projectors. There are three steps: for the $i$-th linear functional $\lambda_i$,
\begin{itemize}
\item [1)] Choose a sub-domain $\Omega_i \subset \Omega$ with the property that $\Omega_i \cap \text{supp} B_i \neq \emptyset$.
\item [2)] Choose a local projector $Q_i$ into $\mathcal{S}(\Omega_i)$ of the form
\begin{equation}
\label{gLocal}
Q_i(f^i):=\sum\limits_{j\in K_i}c_jB_j,~~K_i:=\{j: \Omega_i\cap \text{supp} B_j \neq \emptyset\},
\end{equation}
	where $f^i$ denotes the restriction of $f$ to the sub-domain $\Omega_i$.
\item [3)] Set coefficient $i$ of the global approximation $Q(f)$ to $c_i$, i.e.
	$\lambda_{i}(f)=c_i$.
\end{itemize}
Following the same approach in \cite{examples,LYCHE2008416}, it can be shown that if the local projector satisfies $Q_i(f^i)=f^i,~~\forall f \in \mathcal{S}$,  the quasi-interpolation constructed in this way is also a subdivision space projector when the subdivision space is spanned by linearly independent basis functions.

The local projector $Q_i$ is constructed by local interpolation on the sub-domain $\Omega_i$. Suppose there are $n_i$ non-vanishing basis functions on $\Omega_i$, i.e. $\lvert K_i\rvert=n_i$. We choose $n_i$ interpolation points in $\Omega_i$ such that the local interpolation problem has a unique solution. Then the local interpolation problem is expressed as follows, 
\[\bm{Ac=f},\]
where $\bm{A}$ contains basis values at interpolation points $x_{i,k}\in \Omega_i$,  $\bm{A}=[a_{kj}]$,  $a_{kj}=B_j(x_{i,k})$, $j\in K_i$, $k=1,\cdots,n_i$. The vector $\bm{f}=[f(x_{i,1}),\cdots,f(x_{i,n_i})]^T$ and $\bm{c}$ is the unknown coefficient vector. 
If matrix $\bm A$ is invertible, then the local interpolation problem has a unique solution $\bm {c=A}^{-1}\bm f$ and will reproduce the local subdivision space $\mathcal{S}(\Omega_i)$ . Then we have $\lambda_i(f)=c_i=\bm{\omega f}$, where $\bm\omega$ is some line of $\bm{A}^{-1}$.

For subdivision spaces, the most difficult part is the selection of interpolation points and the calculation of the matrix $\bm{A}$ at extraordinary points, due to the infinite polynomial pieces around extraordinary points. In this paper, we construct $\bm{A}$ based on subdivision schemes rather than the evaluation of subdivision basis functions, facilitating a simple and unified construction framework for the local interpolation problem.

\begin{figure}[htbp]
	\centering
	\subfigure[Vertex mask]{\includegraphics[scale=0.5]{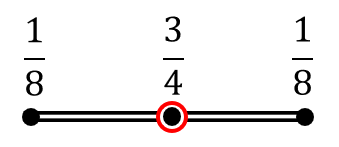}}
	\subfigure[Edge mask]{\includegraphics[scale=0.5]{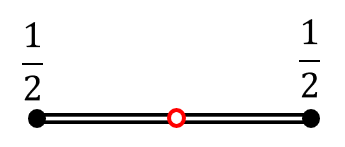}}
	\subfigure[Position mask]{\includegraphics[scale=0.5]{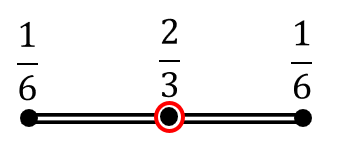}}
	\caption{\label{F:bsplinesubd} The subdivision rules for cubic B-splines, where the circles denote the new vertices.}
\end{figure}

\begin{figure}[htbp]
	\centering
	\includegraphics[scale=0.25]{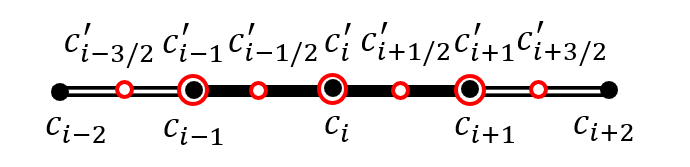}
	\caption{Initial and subdivided control points, where solid line represents local interval $\Omega_i$.}
	\label{F:bsplineAC}
\end{figure}

Let's start with the well-known quasi-interpolation projector of cubic B-splines \cite{examples}. The B-spline of degree $p$ is defined over $p+1$-regular knot vector $\bm{t}=\{t_1,t_2,\cdots,t_{n+p+1}\}$ and denote $N_i^p(t)=N[t_i,t_{i+1},\cdots,t_{i+p+1}](s)$ as the B-spline basis. 
For the $i$-th linear functional $\lambda_i$, it generally chooses $\Omega_i = [t_{i+1}, t_{i+3}]$ as the sub-domain for constructing local projectors.  Because the dimension of the local spline space on $I$ is $5$, one needs to choose five interpolation points to uniquely determine a solution.
For ease of calculation, three knots and two midpoints are chosen, i.e.
 $$x_{i,0}=t_{i+1}, x_{i,1}=(t_{i+1}+t_{i+2})/2, x_{i,2}=t_{i+2}, x_{i,3}=(t_{i+2}+t_{i+3})/2, x_{i,4}=t_{i+3}.$$
The local interpolation problem now has the form $\bm{Ac=f}$ with $\bm{A}=[N_j^3(x_{i,k})]$, $j=i-2,i-1,\cdots,i+2$, $k=0,1,\cdots,4$, $\bm {c}=[c_{i-2},\dots,c_{i+2}]^T$. We show how to construct the matrix $\bm{A}$ in the case of uniform knot vectors based on the B-spline subdivision. The cubic B-spline subdivision rule is shown in Fig.~\ref{F:bsplinesubd}. 
On the local region $I_i$, referring to Fig.~\ref{F:bsplineAC}, we have
\[\bm{Ac=c^\infty:=Lc'=L(Sc)},\]
where $\bm {c'}=[c'_{i-3/2},c'_{i-1},\dots,c'_{i+3/2}]^T$ contains the coefficients after the unknown coefficient vector $\bm{c}$ is subdivided once, thus $\bm{c'=Sc}$, $\bm{S}$ is the subdivision matrix. The limit positions of the coefficients $c'_j,\;j=i-1,i-1/2,\dots,i+1$ are denoted by $\bm {c^\infty}=[c^\infty_{i-1},c^\infty_{i-1/2},\dots,c^\infty_{i+1}]^T$, thus $\bm{c^\infty=Lc'}$, $\bm L$ is position matrix. Meanwhile, $\bm {c^\infty}$ is exactly the value at the interpolation points $x_{i,k}$ of the B-spline determined by the control coefficients in $\bm c$. Then matrix $\bm{A}$ can be constructed by the subdivision matrix $\bm S$ and the position matrix $\bm L$
\[\bm{A=LS}=\frac{1}{6} \begin{bmatrix}
	1&4  &1  &  &  &  & \\
	&1  &4 &1  &  &  & \\
	&  &1  &4  &1  &  & \\
	&  &  &1  &4  &1  & \\
	&  &  &  &1  &4  &1
\end{bmatrix}\cdot \frac{1}{8} \begin{bmatrix}
	4 &4  &  &  & \\
	1 &6  &1  &  & \\
	&4  &4  &  & \\
	&1  &6  &1  & \\
	&  &4  &4  & \\
	&  &1  &6  &1 \\
	&  &  &4  &4
\end{bmatrix}=\frac{1}{48}\begin{bmatrix}
	8 &32  &8  &  & \\
	1&23  &23  &1  & \\
	&8  &32  &8  & \\
	&1  &23  &23  &1 \\
	&  &8  &32 &8
\end{bmatrix}.\]

Since $c_i$ is the third element in $\bm c$, we take the third line of the inverse of $\bm A$ as vector $\bm \omega=[1\;-8\;20\;-8\;1]/6$. Therefore, we have quasi-interpolation projector consisted with \cite{examples}
\[\lambda_{i}f=c_i=\bm{\omega f}=\frac{1}{6}\big(f(t_{i+1})-8f(t_{i+3/2})+20f(t_{i+2})-8f(t_{i+5/2})+f(t_{i+3})\big)\]
where $t_{i+3/2}=(t_{i+1}+t_{i+2})/2$ and $t_{i+5/2}=(t_{i+2}+t_{i+3})/2$.

For general subdivision schemes, the matrix $\bm{A}$ in the local interpolation problem can be constructed as above, i.e., $\bm{A}$ equals the product of the position matrix and subdivision matrix. The matrix $\bm{A}$ constructed in this way has an explicit expression, regardless of the vertex valence to which the basis function corresponds.

\section{Several examples of subdivision space projectors \label{examples}}
In this section, we explain how to construct the quasi-interpolation projectors for the Catmull–Clark, Loop, and modified Loop subdivisions in detail. It is assumed that each face in a control grid contains at most one extraordinary point. This can be easily met by subdividing the control grid once. Here an extraordinary point (\textbf{Ep} for short) means it has other than four patches
adjacent to it in quadrilateral meshes, or it has other than six patches
adjacent to it in triangular meshes. While, a regular point has four patches
adjacent to it in quadrilateral meshes, or it has six patches
adjacent to it in triangular meshes.


\subsection{Catmull–Clark subdivision}
\begin{figure}[htbp]
\centering
 \subfigure[1-ring neighbourhood]{\includegraphics[scale=0.65]{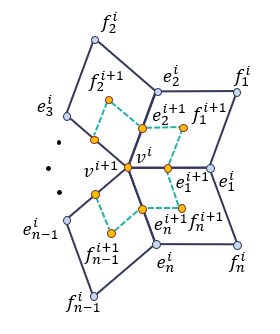}}
  \subfigure[Vertex mask]{\includegraphics[scale=0.65]{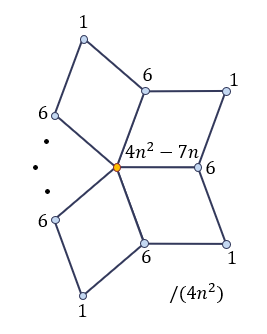}}
\subfigure[Position mask]{\includegraphics[scale=0.65]{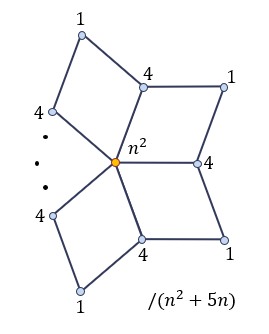}}\\
 \subfigure[Face mask]{\includegraphics[scale=0.65]{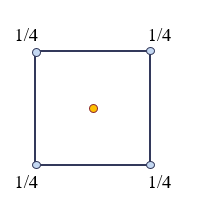}}
 \hspace{0.2cm}
 \subfigure[Edge mask]{\includegraphics[scale=0.65]{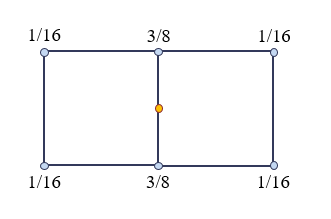}}
\caption{\label{F:ccsubdrules} Catmull–Clark subdivision rules and the limit position mask.}
\end{figure}

Catmull–Clark subdivision \cite{CATMULL1978350} is a quadrilateral-based subdivision scheme. During the subdivision, each face is subdivided into four sub-faces, generating a new face point per face and a new edge point per edge. Let $v^i$ be a vertex in a quadrilateral mesh surrounded by $n$ points $e^i_1,\cdots,e^i_n$ and $n$ faces. Referring to Fig.~\ref{F:ccsubdrules}, the Catmull–Clark subdivision rule is defined as follows: 
\begin{itemize}
	\item Each face point $f_{j}^{i+1}$, $j=1,\cdots,n$ is placed at the centroid of that face,
 \[f_{j}^{i+1}=\frac{v^{i}+e_{j}^{i}+e_{j+1}^{i}+f_{j}^{i}}{4}.\]
	\item Each new edge point $e^{i+1}_j$, $j=1,\cdots,n$ is
computed as an average of four points,
\[e_{j}^{i+1}=\frac{v^{i}+e_{j}^{i}+f_{j-1}^{i+1}+f_{j}^{i+1}}{4}.\]
	\item The vertex $v^i$ is updated to $v^{i+1}$,
 \[v^{i+1}=\frac{n-2}{n}v^i+\frac{1}{n^2}\sum_{j=1}^ne_j^i+\frac{1}{n^2}\sum_{j=1}^nf_j^{i+1}.\]
\end{itemize}

The limit position of the vertex $v^i$ is computed by \cite{Halstead1993} and given by
\[v^\infty=\frac{n^2v^i+4\sum_je_j^i+\sum_jf_j^i}{n(n+5)}.\]

For Catmull-Clark subdivision surfaces, each vertex in the control grid corresponds to a subdivision basis function, and the subdivision basis function is non-vanishing within the 2-ring neighborhood of the associated vertex. Catmull-Clark basis functions are nonnegative, form a partition of unity and have
global linear independency \cite{subdivisionlinearind}.

\begin{figure}[htbp]
	\centering
	\subfigure[Local domain $\Omega_i$ and interpolation points]{\includegraphics[scale=0.25]{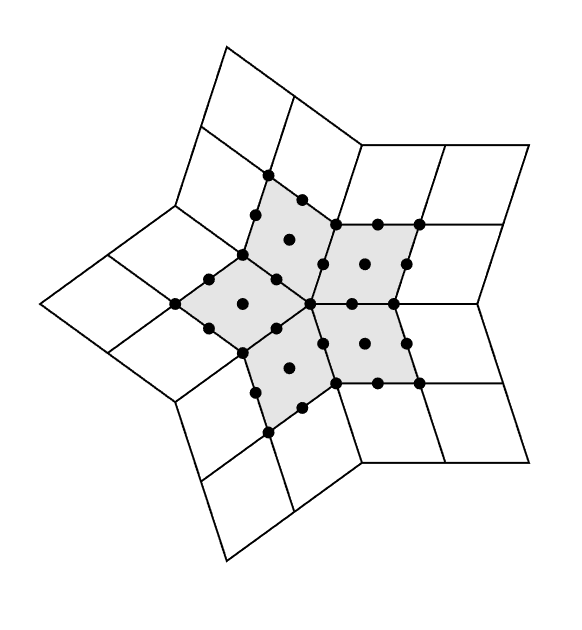}}
	\subfigure[The non-vanishing basis functions on $\Omega_i$]{\includegraphics[scale=0.25]{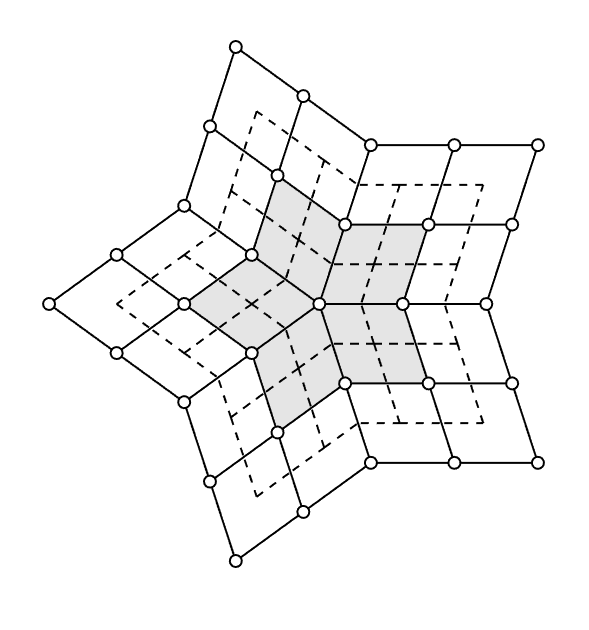}}
	\caption{\label{F:interpcc} The local domain $\Omega_i$ is shaded by gray. (a) The interpolation points around an extraordinary point are marked by black solid circles. (b) The associated vertices of the non-vanishing basis functions on the local domain $\Omega_i$ are marked by white circles. The mesh around $\Omega_i$ after one subdivision is shown by dashed grid.}
\end{figure}

\begin{figure}[htbp]
	\centering
	\subfigure[Extraordinary point case]{\includegraphics[scale=0.18]{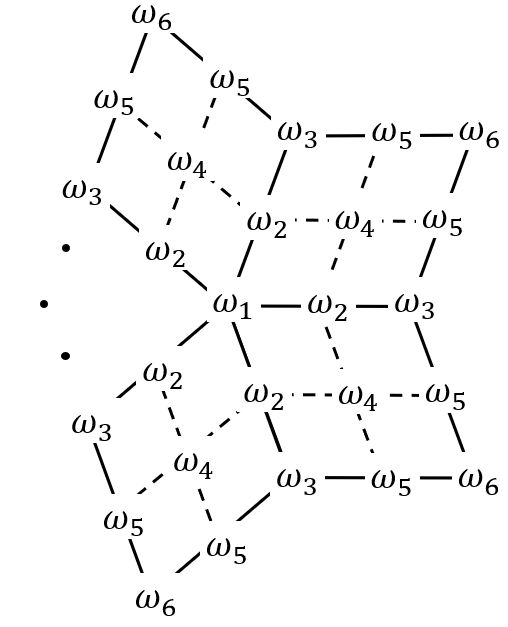}}
	\subfigure[Regular point case]{\includegraphics[scale=0.26]{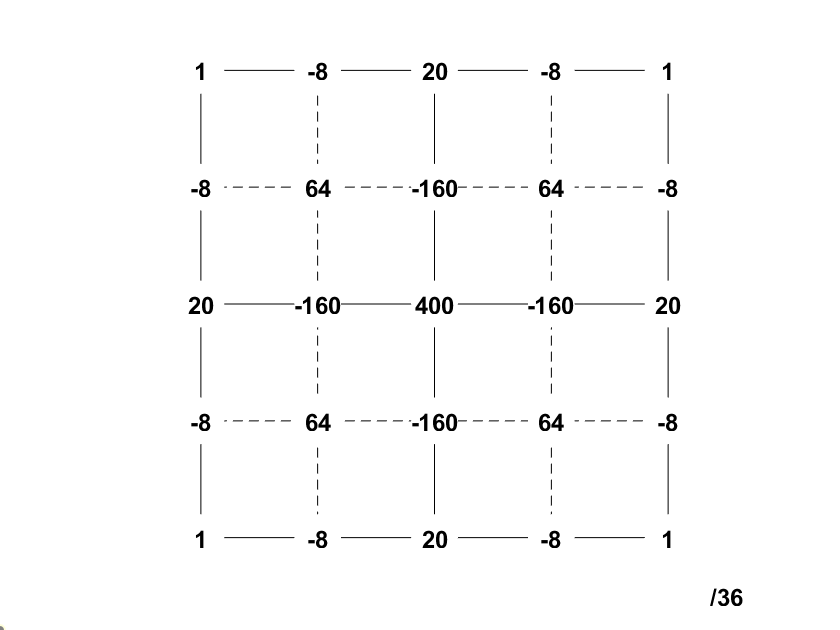}}
	\caption{\label{F:quasicoeffccc} The coefficients of the linear functional $\lambda_i$ defined by \eqref{lamdaK} corresponds to a vertex of valence $n$ for the quasi-interpolation of the Catmull–Clark subdivision.}
\end{figure}

The quasi-interpolation is constructed based on the general recipe in Section \ref{s.generalrecipe}. For a vertex $v$ of valence $n$ ( $n$ patches adjacent to the vertex), the 1-ring neighborhood is chosen as the local region $\Omega_i$ for constructing the local interpolation problem. There are $6n+1$ subdivision basis functions that do not vanish on $\Omega_i$, so $6n+1$ interpolation points,  as shown in Fig.~\ref{F:interpcc}, are selected for the construction of the local interpolation problem $\bm{A}\bm{c}=\bm{f}$ with a unique solution. The interpolation point group consists of all 1-ring neighboring points, edge midpoints, and face centroid points. 
The matrix $\bm A$ is created by multiplying the position matrix $\bm L$ and the subdivision matrix $\bm S$, i.e. $\bm{A=LS}$,
where the number of rows in $\bm S$ and the number of columns in matrix $\bm L$ are equal to $12n+1$. The $\bm L$ and $\bm S$ as given in detail in Appendix A. Then, we obtain
\[\lambda_{i}(f)=\bm{\omega f},\]
where $\bm{\omega}$ is the $k$-th row of ${\bm A}^{-1}$ and $k$ is the index  number of the vertex $v$  in the vector $\bm{c}$. There are only 6 different values in $\bm{\omega}$, which are denoted by $\omega_i,\;i=1,2,\dots,6$, and
\begin{equation}
	\left\{
	\begin{aligned}
		\omega_{6} &=\frac{18(5n-2)}{n(739n^{2}-2717n+1960)} \\
		\omega_{3} &=\frac{18(91n-4)}{n(739n^{2}-2717n+1960)} \\
		\omega_{5} &= -8\omega_{6},\quad \omega_{2} = -8\omega_{3}, \\
		\omega_{4}&=64\omega_{6},\quad \omega_{1} = 1-n(49\omega_{6}-7\omega_{3}).
	\end{aligned}
	\right.
\end{equation}
The quasi-interpolation scheme for the Catmull-Clark subdivision is concluded in Fig.~\ref{F:quasicoeffccc}.

Especially, if $n=4$, i.e $v$ is a regular point, we have
\[\omega_1=\frac{100}{9},\quad\omega_2=-\frac{40}{9},\quad\omega_3=\frac{5}{9},\quad
\omega_4=\frac{16}{9},\quad\omega_5=-\frac{2}{9},\quad\omega_6=\frac{1}{36}.\]
Referring to Fig.~\ref{F:quasicoeffccc}(b), the linear functional $\lambda_i(f)$ for $n=4$ is the same as the tensor product of the linear functional deduced from uniform cubic B-splines (see Section \ref{s.generalrecipe}). 


\subsection{Loop subdivision \label{s.loopQI}}
\begin{figure}[htbp]
	\centering
        \subfigure[1-ring neighbourhood]{\includegraphics[scale=0.3]{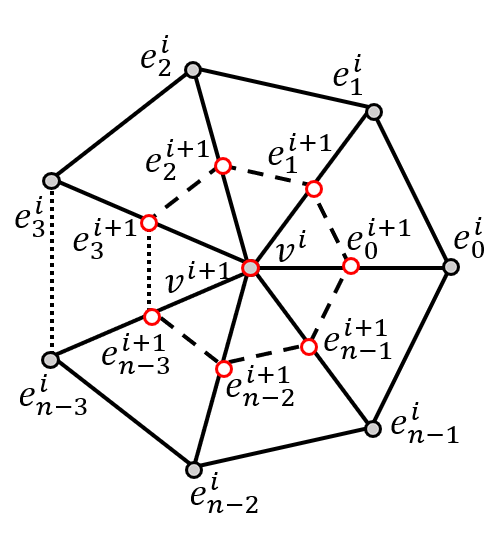}}
	\subfigure[Vertex mask]{\includegraphics[scale=0.3]{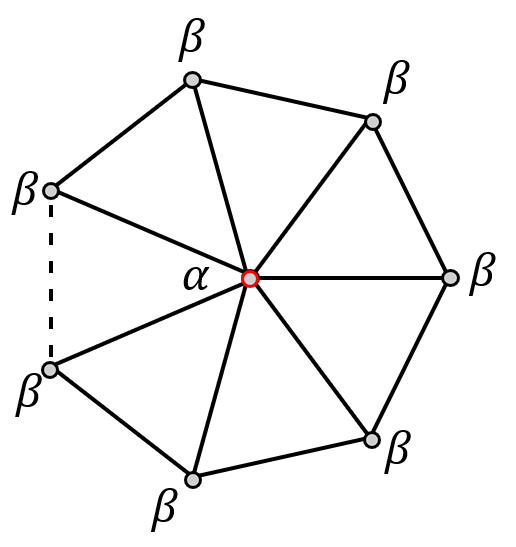}}
	\subfigure[Edge mask]{\includegraphics[scale=0.3]{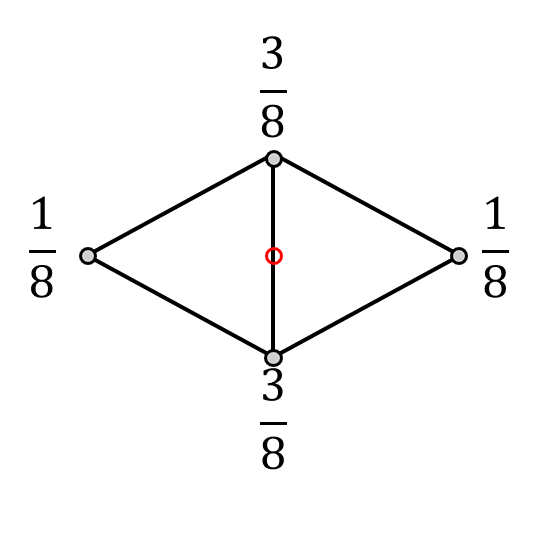}}
	\subfigure[Position mask]{\includegraphics[scale=0.3]{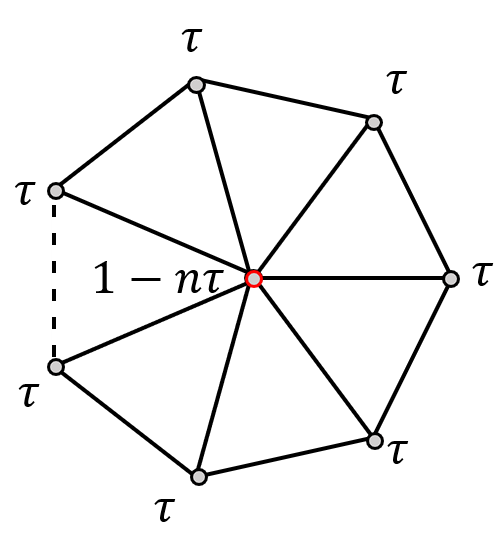}}
	\caption{\label{F:loopsubd} Loop subdivision rules and the limit position mask.}
\end{figure}

Loop subdivision \cite{Loop} is a triangle-based subdivision scheme. During the subdivision, each triangle is subdivided into four sub-triangles, generating a new point per edge. Referring to Fig.~\ref{F:loopsubd}, let $v^i$ be a valence $n$ vertex and $e_j^i$, $j=0,1, \cdots, n-1$ be the $n$ neighboring vertices. Each edge $v^ie^i_j$ adds a new edge point $e^{i+1}_j$ and the vertex $v^i$ is replaced by a new vertex $v^{i+1}$, where
\begin{eqnarray}
   v^{i+1} &=&\alpha v^{i}+\beta \sum_{j=0}^{n-1} e^i_j,~~\alpha = \dfrac{3}{8}-\left(\dfrac{3}{8}+\dfrac{1}{4}\cos\dfrac{2\pi}{n}\right)^2,~~\beta = \dfrac{1-\alpha}{n},\\
    e^{i+1}_j&=&\frac{1}{8}(e_{j-1}^i+e_{j+1}^i) + \frac{3}{8}(e_j^i+v^i).
\end{eqnarray}
According to the position mask given in \cite{MA2004525}, the limit position of a vertex $v^i$ is computed by
\[v^{\infty}=(1-n\tau)v^i+\tau\sum_{j=0}^{n-1}e_{j}^i,\text{ where }\tau=\left(\frac{3}{8\beta}+n\right)^{-1}.\]  

\begin{figure}[htbp]
	\centering
	\subfigure[Interpolation points]{\includegraphics[scale=0.25]{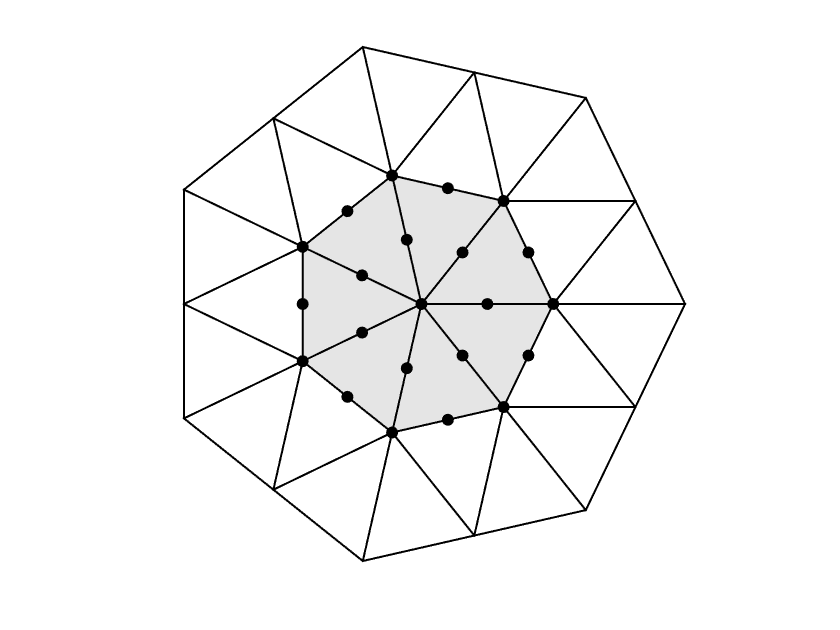}}
	\subfigure[Non-vanishing basis functions]{\includegraphics[scale=0.25]{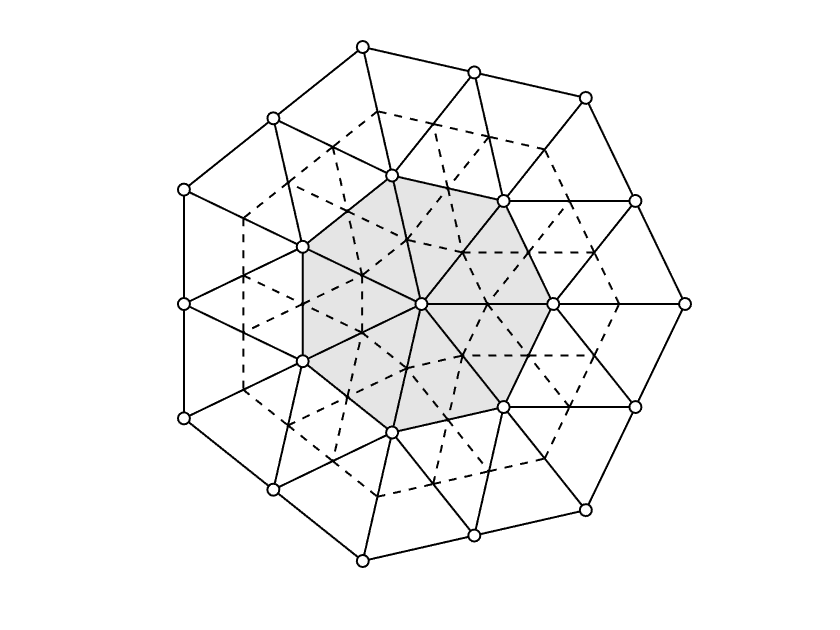}}
	\caption{\label{F:interploop} The local domain $\Omega_i$ is shaded by gray. (a) The interpolation points around an extraordinary point are marked by black solid circles. (b) The associated vertices of the non-vanishing basis functions on the local domain $\Omega_i$ are marked by white circles. The mesh around $\Omega_i$ after one subdivision is shown by dashed grid. }
\end{figure}

For Loop subdivision surfaces, each vertex in the control grid corresponds to a subdivision basis function, and the subdivision basis function is non-vanishing within the 2-ring neighborhood of the associated vertex. Loop basis functions are nonnegative, form a partition of unity and have global linear independency \cite{subdivisionlinearind}.

\begin{figure}[htbp]
	\centering
	\subfigure[Extraordinary point case]{\includegraphics[scale=0.20]{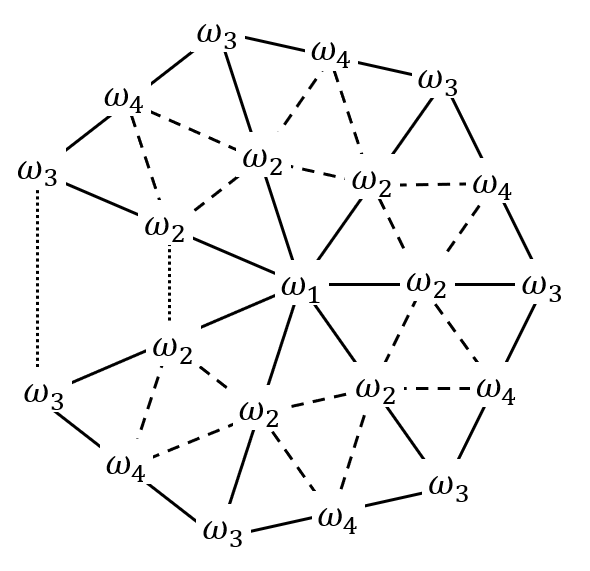}}
	\subfigure[Regular point case]{\includegraphics[scale=0.26]{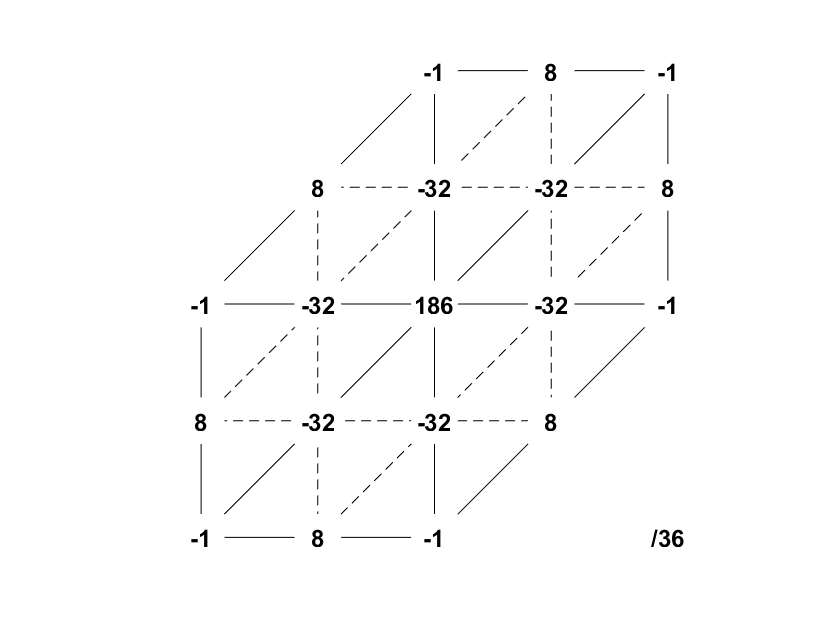}}
	\caption{\label{F:quasicoeffloop} The coefficients of the linear functional $\lambda_i$ defined by \eqref{lamdaK} corresponds to a vertex of valence $n$ for the quasi-interpolation of the Loop subdivision.}
\end{figure}

The quasi-interpolation is constructed based on the general recipe in Section \ref{s.generalrecipe}. For a vertex $v$ of valence $n$ ( $n$ patches adjacent to the vertex), the 1-ring neighborhood is chosen as the local region $\Omega_i$ for constructing the local interpolation problem. There are $3n+1$ subdivision basis functions that do not vanish on $\Omega_i$, so $3n+1$ interpolation points,  as shown in Fig.~\ref{F:interploop}, are selected for the construction of the local interpolation problem $\bm{A}\bm{c}=\bm{f}$ with a unique solution. The interpolation point group consists of all 1-ring neighboring points, edge midpoints, and face centroid points. 
The matrix $\bm A$ is created by multiplying the position matrix $\bm L$ and the subdivision matrix $\bm S$, i.e. $\bm{A=LS}$,
where the number of rows in $\bm S$ and the number of columns in matrix $\bm L$ are equal to $6n+1$. The $\bm L$ and $\bm S$ is given in detail in Appendix A. Then, we obtain
\[\lambda_{i}(f)=\bm{\omega f},\]
where $\bm{\omega}$ is the $k$-th row of ${\bm A}^{-1}$ and $k$ is the index  number of the vertex $v$  in the vector $\bm{c}$. There are only 4 different values in $\bm{\omega}$, which are denoted by $\omega_i,\;i=1,\dots,4$, and
\begin{equation}
	\omega_{3}=\frac{-48\beta}{128\alpha^{2}+200\alpha-67},\quad \omega_{2}=32\omega_{3},\quad \omega_{4}=-8\omega_{3},\quad \omega_{1}=1-25n\omega_{3}  	
\end{equation}
The quasi-interpolation scheme for the Loop subdivision is concluded in Fig.~\ref{F:quasicoeffloop}.

Especially, if $n=6$, the regular case, we have
\[\omega_1=\frac{31}{6},\quad\omega_2=-\frac{8}{9},\quad\omega_3=-\frac{1}{36},\quad
\omega_4=\frac{2}{9}.\]
Referring to Fig.~\ref{F:quasicoeffloop}(b), the linear functional $\lambda_i(f)$ for $n=6$ is the same as the quasi-interpolation of the quartic $C^2$-continuous box splines discussed in \cite{LYCHE2008416}.

\subsection{Modified Loop subdivision\label{s.loopGQI}}
A modified version of the Loop subdivision was proposed in \cite{KANGmodifiedloop} to achieve the optimal approximation rate in isogeometric analysis. The subdivision rule at extraordinary points is controlled by a parameter $\lambda$. Referring to Fig.~\ref{F:loopsubd}, let $v^i$ be a valence $n$ vertex and $e_j^i$, $j=0,1, \cdots, n-1$ be the $n$ neighboring vertices. Each edge $v^ie^i_j$ adds a new edge point $e^{i+1}_j$ and the vertex $v^i$ is replaced by a new vertex $v^{i+1}$, where
\begin{equation}
\label{E.modifiedloop}
     v^{i+1} =\alpha v^{i}+\beta \sum_{j=0}^{n-1} e^i_j,~~
    e^{i+1}_k=\gamma v^{i} + \sum_{j=0}^{n-1} \gamma_je^i_{j+k},~~\alpha+n\beta=\gamma+\sum_{j=0}^{n-1} \gamma_j=1,
\end{equation}
where $\beta$ and $\gamma_j$ are functions of the parameter $\lambda$ and given in \cite{KANGmodifiedloop}. Note that only extraordinary points and the edge points adjacent to extraordinary points use a different subdivision rule from the Loop subdivision, and the remaining points are updated by the Loop subdivision.

\begin{figure}[htbp]
\centering
\subfigure[Edge mask]{\includegraphics[scale=0.35]{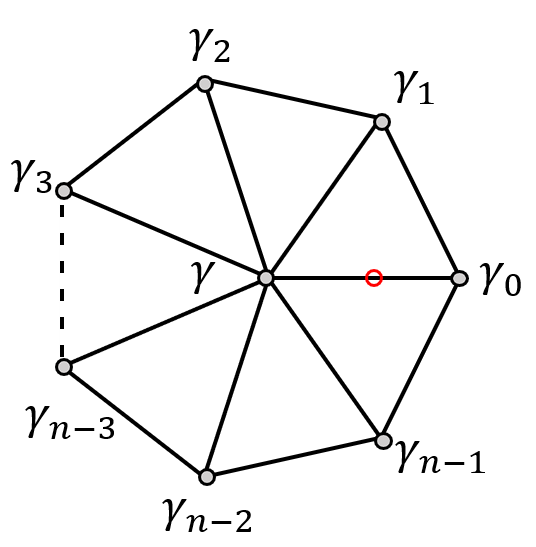}} \hspace{0.2cm}
\subfigure[Position mask for extraordinary points]{\includegraphics[scale=0.35]{loopposition.png}}\hspace{0.2cm}\\
\subfigure[Position mask for regular points]{\includegraphics[scale=0.4]{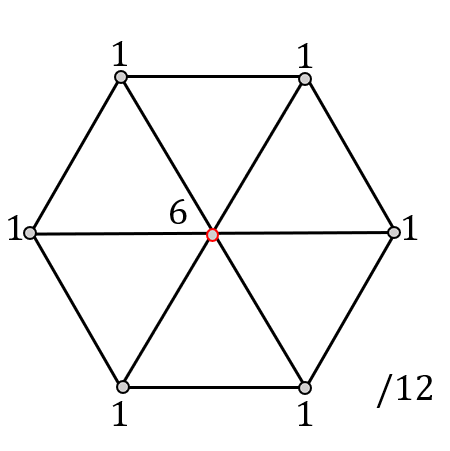}}
\subfigure[Position mask for regular points adjacent to Eps]{\includegraphics[scale=0.4]{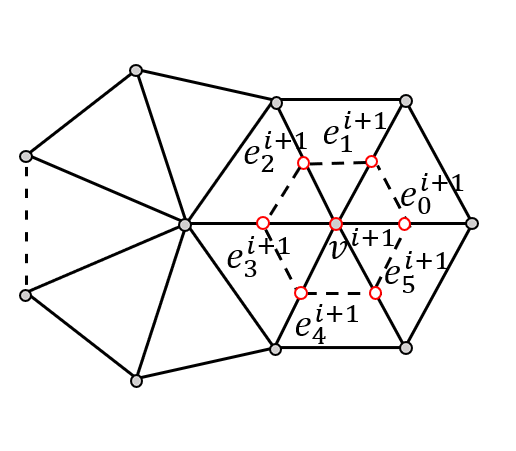}
\includegraphics[scale=0.4]{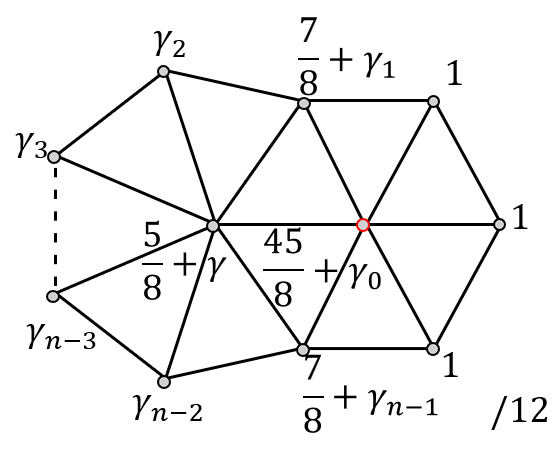}}
\caption{Subdivision rules and position masks for modified Loop subdivision.}
	\label{F:Modifiedloopsubd}
\end{figure}

For the modified Loop subdivision, the limit position of a control point can be derived from the eigenvalues of the subdivision matrix. For an extraordinary point $v^i$ (valence $n \neq 6$), the limit position is computed by 
\[v^{\infty}=(1-n\tau)v^i+\tau\sum_{j=0}^{n-1}e^i_{j},\text{ where }\tau=\left(\frac{\gamma}{\beta}+n\right)^{-1}.\] 
This formula is equivalent to the result given in \cite{Boundedcurvature}.  
For a regular point $v^i$ (valence $n=6$), if there are no extraordinary points in its 1-ring neighborhood, the limit position is computed by  
\[v^{\infty}=\frac{1}{2}v^i+\frac{1}{12}\sum_{j=0}^{5}e^i_{j},\]
as shown in Fig.~\ref{F:Modifiedloopsubd}(c). Otherwise, a local subdivision is performed before calculation, as shown in Fig.~\ref{F:Modifiedloopsubd}(d), then the limit position is computed based on the new points generated by one step of subdivision, that is
\[v^\infty=\frac{1}{2}v^{i+1}+\frac{1}{12}\sum_{j=0}^{5}e^{i+1}_{j}=\sum_{j=0}^{n+2}\theta_j e^i_j+\theta_{n+3}v^i,\]
where $e^{i+1}_3$  is computed by the edge mask Fig.~\ref{F:Modifiedloopsubd}(a) and other points $e^{i+1}_{j}$, $v^{i+1}$ are updated based on Loop subdivision. By substituting the expressions of $e^{i+1}_{j}$ and $v^{i+1}$, the limit position becomes a linear combination of $e^i_j,\;j=0,\dots,n+2$ and $v^i$, where $n$ is the valence of the adjacent extraordinary point. The coefficients $\theta_j$ are given in Fig.~\ref{F:Modifiedloopsubd}(d).

In a triangular mesh, each vertex corresponds to a modified Loop subdivision basis function. Compared to the Loop subdivision, the supports of some modified Loop subdivision basis functions are enlarged and the number of non-vanishing basis functions on an element is also increased because the edge point rule \eqref{E.modifiedloop} involves all the 1-ring neighbors of an extraordinary point. However, the basis functions are still linearly independent and have local support.

The local domain and the interpolation points are selected in the same way as the Loop subdivision. The subdivision matrix $\bm{S}$ and the position matrix $\bm{L}$ are given in Appendix A in detail. Then we have $\lambda_{i}(f)=\bm{\omega f}$ by solving $\bm{Ac =f}$. 
Similar to the Loop subdivision, we have
\[\omega_{2}=32\omega_{3},\quad \omega_{4}=-8\omega_{3},\quad \omega_{1}=1-25n\omega_{3}. \] 
The different one $\omega_3$ is obtained by 
\begin{equation}
    \omega_{3}=\frac{-96\beta}{P_{1}(\gamma)\alpha^{2}+P_{2}(\gamma)\alpha+P_{3}(\gamma)},
\end{equation}
where
\[\left\{\begin{aligned}
	&P_{1}(x)=256x+160 \\
	&P_{2}(x)=-512x^{2}+1608x-131 \\
	&P_{3}(x)=256x^{3}-1768x^{2}+347x-29
\end{aligned}\right. .\]

\section{Treatment of boundary points and regular points adjacent to Eps \label{Sec4}}
Note that the framework for constructing linear functional $\lambda_i(f)$ proposed in Section \ref{examples} applies to extraordinary points and those regular points that are not in the 1-ring neighborhood (or not in the 2-ring neighborhood for the modified Loop subdivision) of extraordinary points. In this section, we discuss how to design linear functional for boundary points and the regular points that are not discussed in Section \ref{examples}.

\subsection{Regular points adjacent to Eps \label{s.adjaeps}}

When we use the method in Section \ref{examples} to obtain a linear functional $\lambda_i(f)$ corresponding to an extraordinary point, we can also get the linear functional corresponding to any 2-ring neighborhood point of that point, which corresponds to a row of the inverse of the matrix $\bm A$. Based on this idea we give the linear functional for the regular points that are adjacent to Eps.

Suppose $v_r$ is a regular point in the 1-ring neighborhood of Ep denoted by $v^r_{ep}$. Intuitively, one can choose the interpolation points of $v^r_{ep}$ as the interpolation points of $v_r$ and solve the coefficients $\bm{\omega}$ associated with $v_r$ from the local interpolation problem related to $v^r_{ep}$. However, the expression of $\bm{\omega}$ is complicated and varies with the valence of $v^r_{ep}$. Therefore, an economical way is to search for a target point within the 1-ring neighborhood (some subdivisions may require a 2-ring neighborhood) of $v_r$: it is regular and has no Eps in its 1-ring neighborhood. In this way, $v_r$ and its target point share the same local interpolation problem $\bm{Ac=f}$, but correspond to different rows of the matrix $\bm{A}$. The details of solution are given in Appendix B.

\begin{figure}[htbp]
	\centering
	\subfigure[Catmull–Clark subdivison]{\includegraphics[scale=0.36]{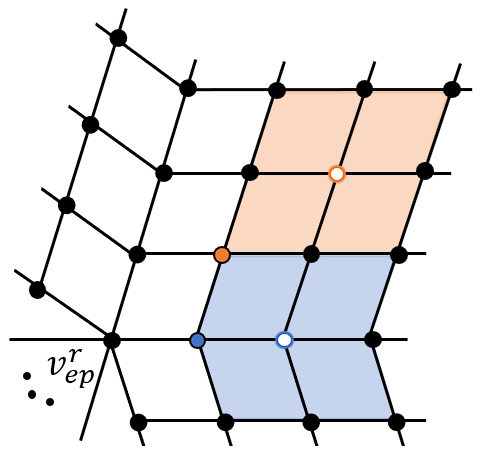}}
 \hspace{0.2cm}
 \subfigure[Loop subdivison]{\includegraphics[scale=0.4]{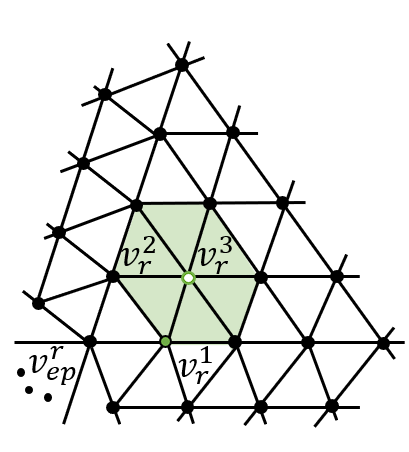}}\\
  \subfigure[Modified Loop subdivision (1-ring neighborhood of $v^r_{ep}$)]{\includegraphics[scale=0.4]{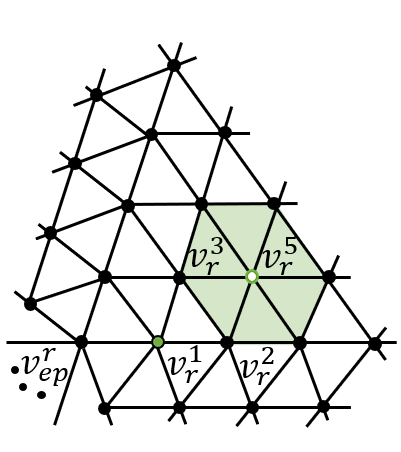}}
  \hspace{0.2cm}
  \subfigure[Modified Loop subdivision (2-ring neighborhood of $v^r_{ep}$)]{\includegraphics[scale=0.4]{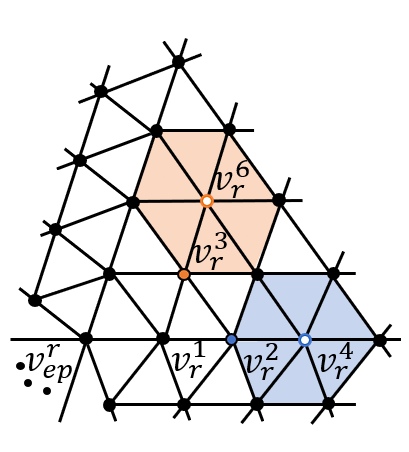}}
\caption{\label{Fig_7} The target point of a regular point (marked with a colored solid circle) is marked with a circle of the same color. The interpolation points are chosen in the shaded domain, as discussed in Section \ref{examples}. Regular points and their target points share the same interpolation points.}
\end{figure}

Note that the target point of $v_r$ may not be unique. We generally choose the one with symmetry, as shown in Fig.~\ref{Fig_7}. This strategy often gives
smaller norm of the quasi-interpolation operator $Q$ and more overlapping interpolation points between adjacent points, resulting in smaller fitting errors.  


For the Catmull–Clark subdivision, the 1-ring neighborhood points of an Ep are classified into two categories: those that share an edge with the Ep, and the remaining ones. The strategy for choosing the target points for these two categories is shown in Fig.~\ref{Fig_7}(a). 
For the Loop subdivision, for the point $v^r_1$, we choose the element that contains the edge $v^r_1v_{ep}^r$ and the element is denoted by $v^r_{ep}v_r^1v_r^2$. We then set the point that is symmetrical with $v_{ep}^r$ relative to the edge $v_r^1v_r^2$ as the target point of $v_r^1$, as the one $v_r^3$ shown in Fig.~\ref{Fig_7}(b). For the modified Loop subdivision, both the 1-ring and 2-ring neighborhood points require the target points. Moreover, there are no Eps in the 2-ring neighborhood of the target points. As shown in Fig.~\ref{Fig_7}(c-d), the target point of $v_r^1$ is denoted by $v_5^r$, and the target point of the 2-ring neighborhood point $v_r^2(v_r^3)$ is chosen as $v_r^4(v_r^6)$.

\begin{remark}
When the target points can not be found in the above way, the adjacent Eps will be served as the target points.
\end{remark}


\subsection{Treatment for boundary points }

For the subdivision of meshes with boundaries, the boundary treatment proposed in \cite{KANGmodifiedloop} is considered here. An element is called a boundary element if it contains at least one boundary vertex, otherwise, it is called an interior element. The subdivision surface is defined over interior elements. This can be understood in the analogy of control meshes for B-spline surfaces. The elements served for defining subdivision surfaces are called active elements. There is no difference between the evaluation of the boundary patches and interior regular patches. The subdivision basis functions that do not vanish on the boundary are all retained, ensuring the complete subdivision space on boundary patches. 

\begin{figure}[htbp]
	\centering
	\subfigure[Corner case (C-C)]{\includegraphics[scale=0.30]{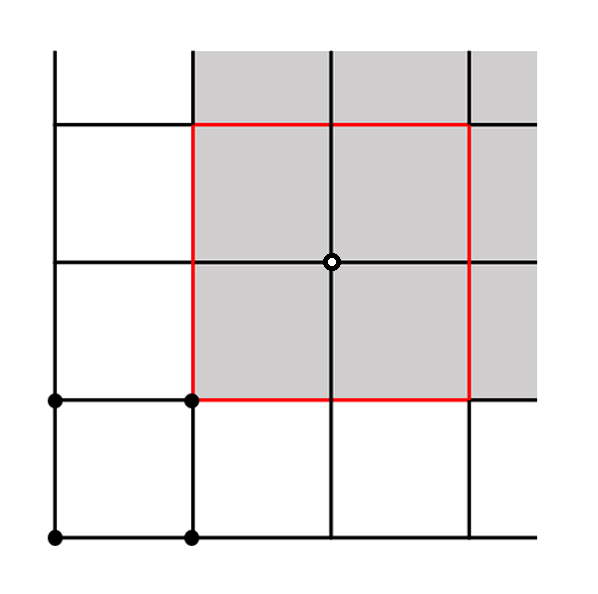}}
	\subfigure[Edge case (C-C)]{\includegraphics[scale=0.30]{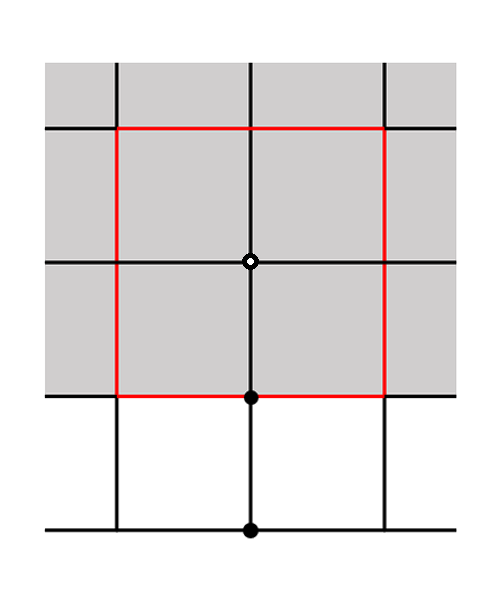}}
	\subfigure[Corner case (Loop)]{\includegraphics[scale=0.33]{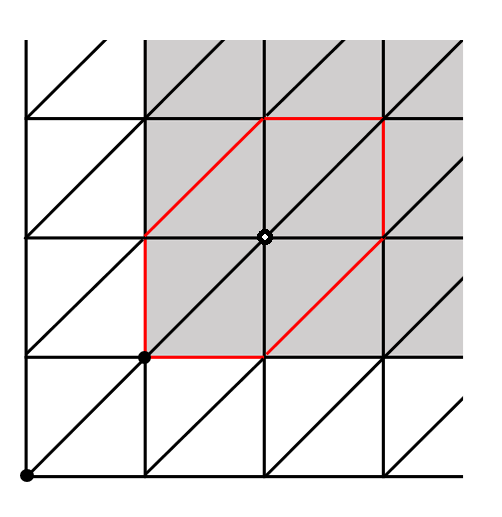}}
	\subfigure[Edge case (Loop)]{\includegraphics[scale=0.33]{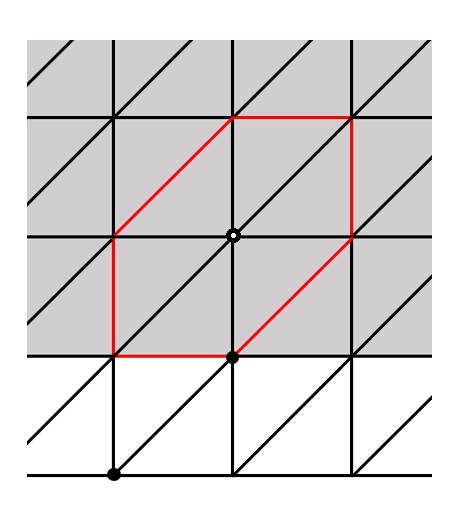}}
	\caption{\label{Fig_9} Illustration of the target points for the boundary points which are marked with black solid circles. The active elements are shaded grey. The interpolation points are selected from an area enclosed by colored polygons. The points on the outermost layer are the mesh boundary points. }
\end{figure}

Based on this boundary processing, there are two layers of points in the mesh that need to be reprocessed: the boundary points of the mesh and the boundary points of the subdivision surfaces (i.e. the boundary points of the active elements). The approach in Section \ref{examples} is not suitable for these two types of boundary points, as it results in some of the interpolated points that are not in the active domain. Therefore, we adopt the same idea in Section \ref{s.adjaeps} to deal with such boundary points.
Fig.\ref{Fig_9} illustrates a way to select target points for mesh and surface boundary points.

\begin{figure}[htbp]
	\centering
	\subfigure[Triangular mesh]{\includegraphics[scale=0.20]{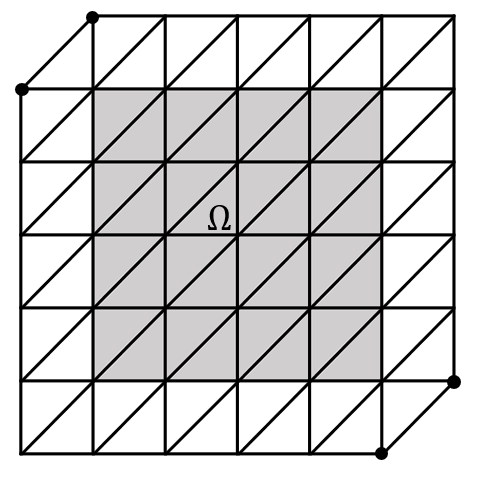}}
	\subfigure[Interpolation points]{\includegraphics[scale=0.18]{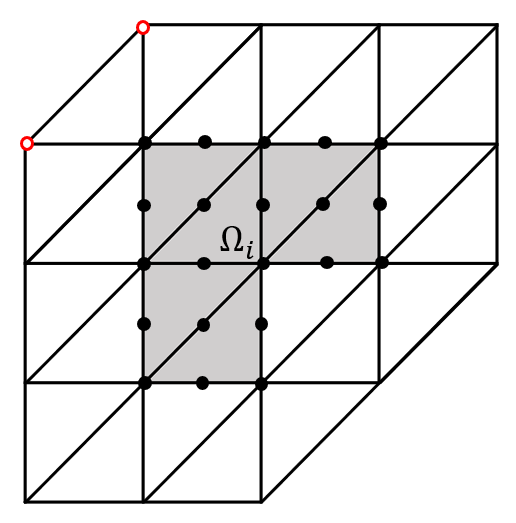}}
	\caption{\label{Fig_10}(a) The target points corresponding to the boundary points marked with solid circles do not exist. (b) For the boundary point (marked with a red circle), the local region $\Omega_i$ is colored in gray and 21 interpolation points are marked with solid circles.}
\end{figure}

Fig.\ref{Fig_10}(a) shows an example that the boundary points can not find the target points. In this case, we can select the local region $\Omega_i$ as shown in Fig.\ref{Fig_10}(b), so that the number of interpolation points is the same as the number of non-vanishing basis functions on the region. Then we construct the matrix $\bm A$ and obtain the linear functionals $\lambda_i(f)$. The details are given in Appendix B.

\section{Numerical Experiment \label{nums}}
\begin{figure}[htbp]
	\centering
\subfigure[Valance 3]{\includegraphics[scale=0.17]{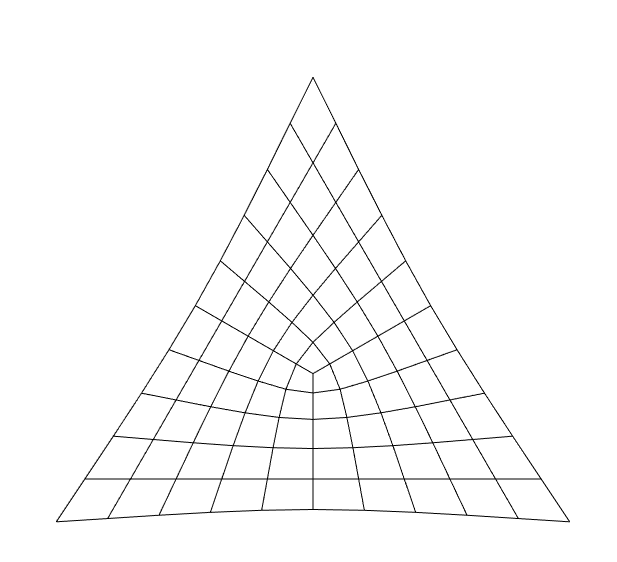}}
\subfigure[Valance 5]{\includegraphics[scale=0.17]{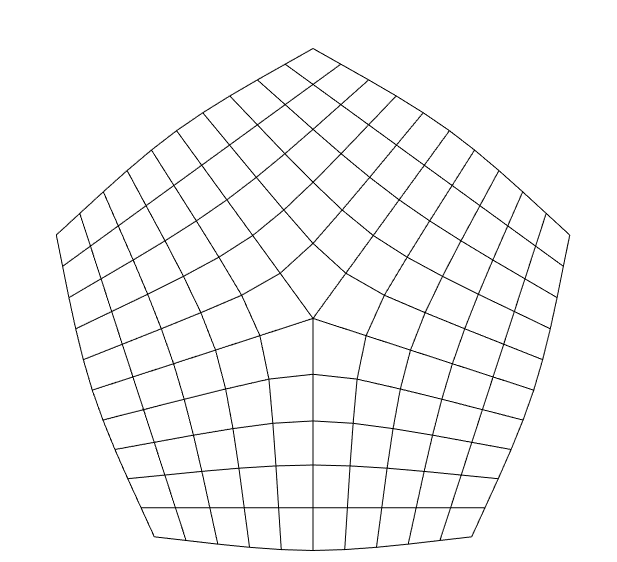}}
\subfigure[Valance 6]{\includegraphics[scale=0.17]{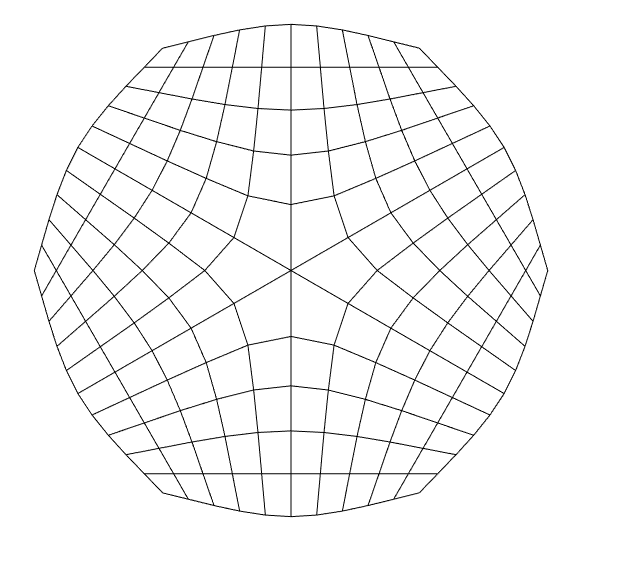}}
\subfigure[Valance 3]{\includegraphics[scale=0.17]{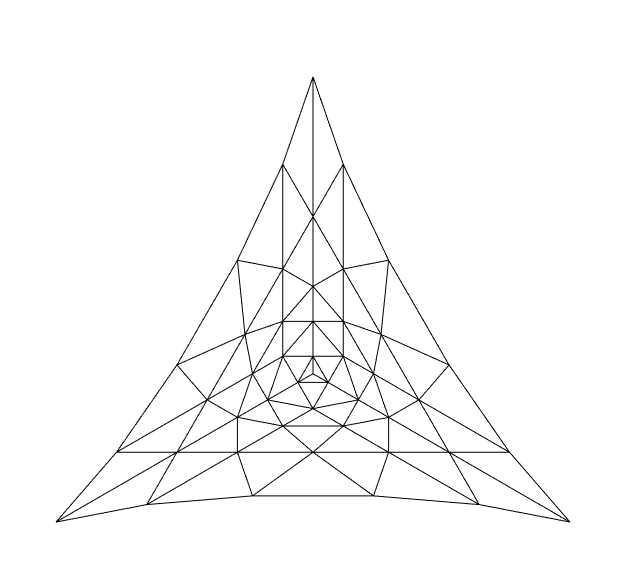}}
\subfigure[Valance 4]{\includegraphics[scale=0.17]{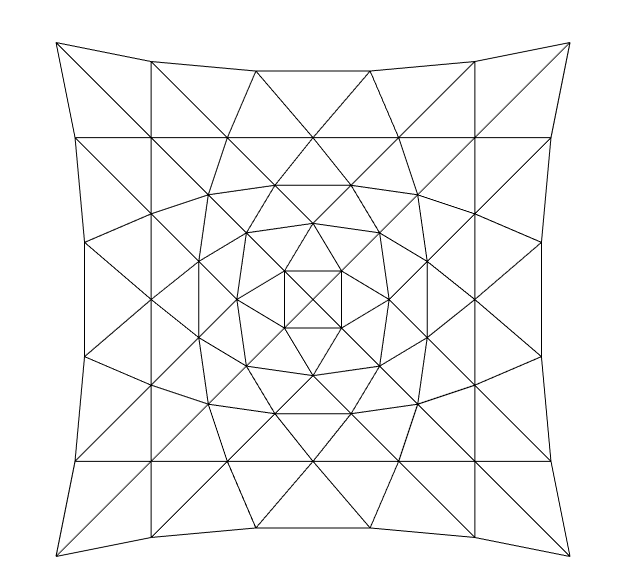}}
\subfigure[Valance 5]{\includegraphics[scale=0.17]{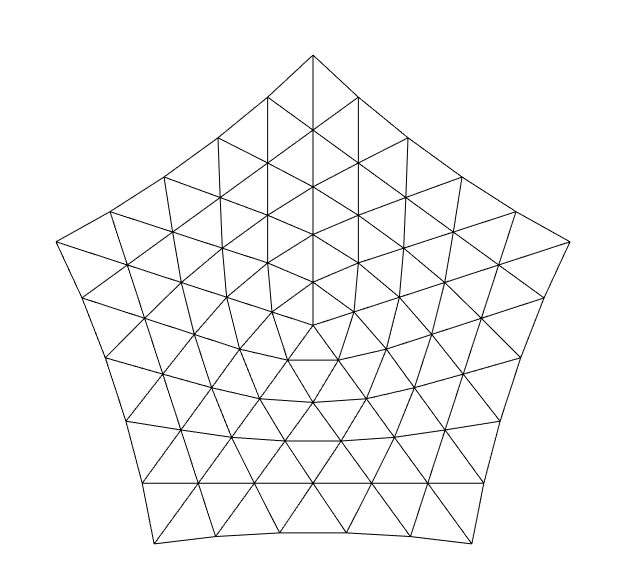}}
\subfigure[Valance 7]{\includegraphics[scale=0.17]{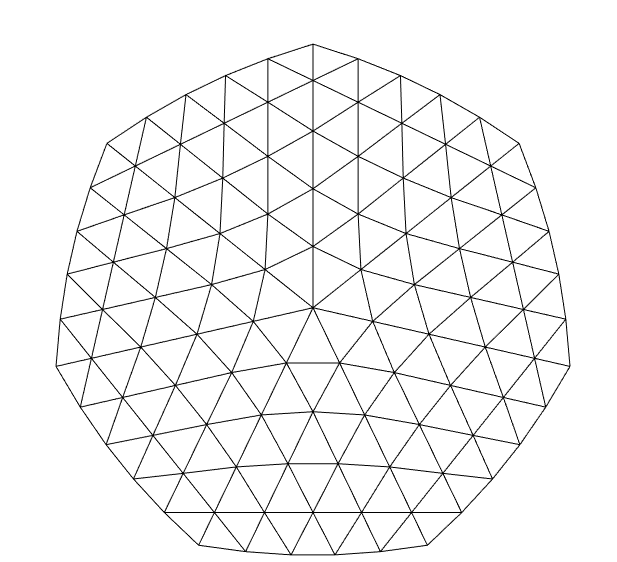}}
\subfigure[Valance 8]{\includegraphics[scale=0.17]{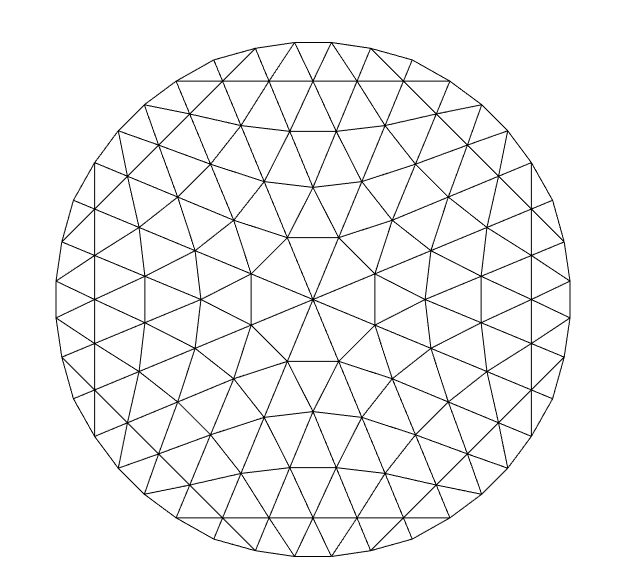}}
	\caption{\label{Fig_16} Several initial meshes with different valences.}
\end{figure}

In this section, we demonstrate the performance of the proposed quasi-interpolation for subdivision schemes. The test function is defined by

\[f(x,y)=e^{-6(x^2+y^2)}, ~~(x,y)\in \Omega, \]
where $\Omega$ represents the subdivision surfaces defined by the initial control meshes shown in  Fig.~\ref{Fig_16}. 
The approximation error is defined by the relative $L_2$ error and relative $L_\infty$ error, i.e.
\begin{equation}
    E_2=\|Q(f)-f\|_2/\|f\|_2,\quad E_\infty=\|Q(f)-f\|_\infty/\|f\|_\infty,
\end{equation}
where 
$\|f\|_2=(\int_\Omega f^2d\bm x)^{\frac{1}{2}},\;\|f\|_\infty=\sup_{\bm x\in\Omega} \lvert f\rvert.$

\begin{figure}[htbp]
	\centering
\subfigure[$L_2$-norm]{\includegraphics[scale=0.26]{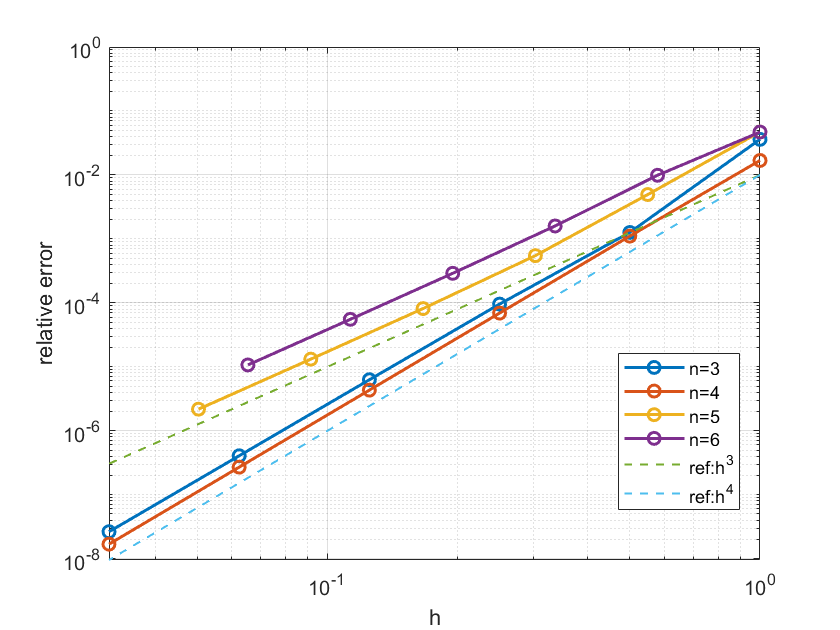}}
\subfigure[$L_\infty$-norm]{\includegraphics[scale=0.26]{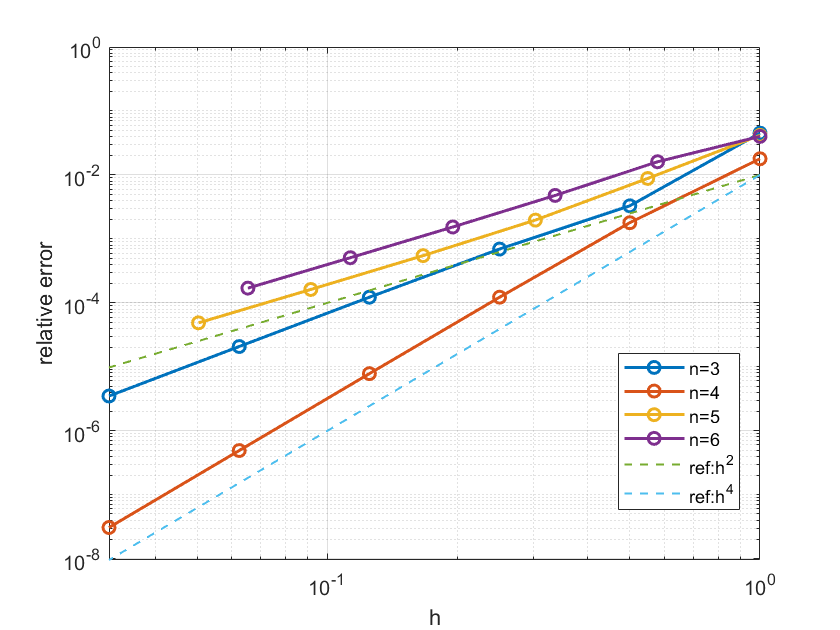}}
	\caption{\label{Fig_19} Relative $L_2$-norm errors and relative $L_\infty$-norm errors from the quasi-interpolation projectors for the Catmull–Clark subdivision on the subdivision surfaces with different valances. The symbol $n$ denotes the valence of the Ep in the mesh shown in Fig.~\ref{Fig_16}(a-c).}
\end{figure}

\begin{figure}[htbp]
	\centering
\subfigure[$L_2$ norm for $n\leq6$]{\includegraphics[scale=0.26]{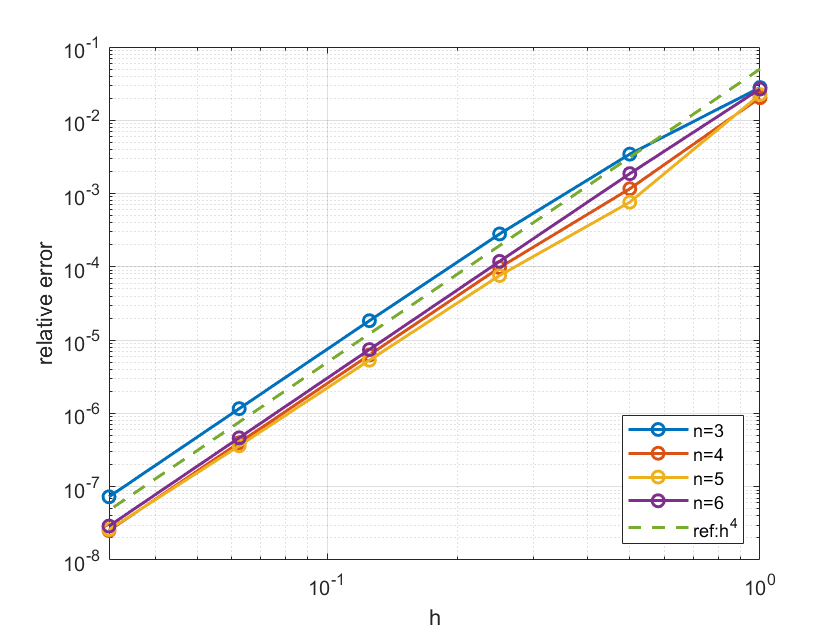}}
\subfigure[$L_\infty$ norm for $n\leq6$]{\includegraphics[scale=0.26]{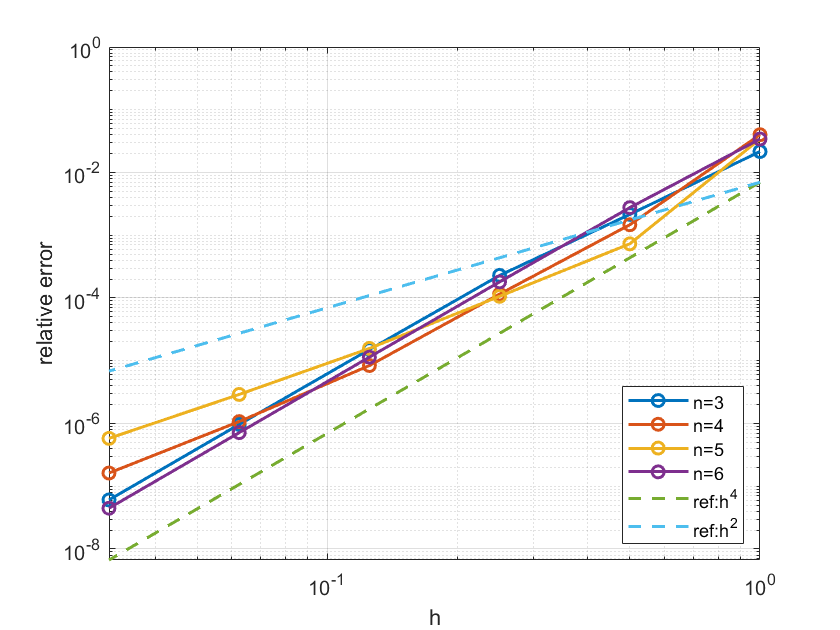}}
\subfigure[$L_2$ norm for $n=7,8$]{\includegraphics[scale=0.26]{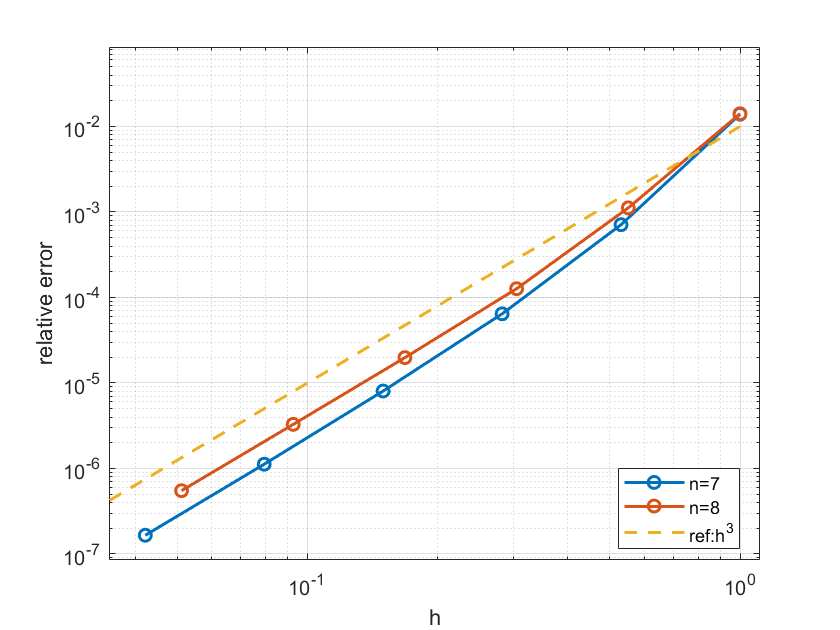}}
\subfigure[$L_\infty$ norm for $n=7,8$]{\includegraphics[scale=0.26]{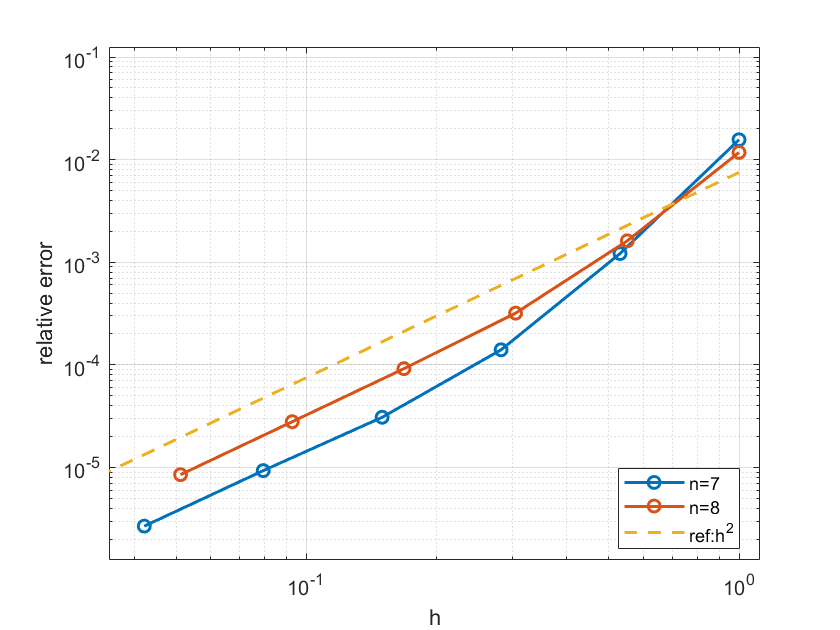}}
	\caption{\label{Fig_21} Relative $L_2$-norm errors and relative $L_\infty$-norm errors from the quasi-interpolation projectors for the Loop subdivision on the subdivision surfaces with different valances. The symbol $n$ denotes the valence of the Ep in the mesh shown in Fig.~\ref{Fig_16}(d-h).}
\end{figure}

Fig.\ref{Fig_19} and Fig.\ref{Fig_21} show the approximation plots of the quasi-interpolation projectors of the Catmull–Clark and Loop subdivision on the subdivision surfaces defined by the initial meshes given in Fig.~\ref{Fig_16}. 
On the mesh without Eps, both the proposed Catmull–Clark and Loop subdivision projectors achieve the optimal approximation order, as now the corresponding subdivision surfaces degenerate to uniform bicubic B-splines and quartic $C^2$-continuous box splines, respectively. On the meshes with one Ep, both the Catmull–Clark and Loop subdivision projectors have the approximation order of 2 under the $L_\infty$-norm and 3 under the $L_2$-norm.

We further discuss the approximation rates of the subdivision space projectors under the $L_2$-norm. In \cite{2001Approximation}, it was proved that for any smooth function $f$, it can be approximated by the Loop subdivision space with an approximation order of $\lambda_{max}^{(3-\epsilon)k}$, where $k$ is the multiresolution level (subdivision level) of the approximation, $\epsilon$ is an arbitrary positive number, and $\lambda_{max}=\max \{\lambda, \frac{1}{2}\}$ with $\lambda$ be the subdominant eigenvalue of the subdivision matrix for valence $n$. This means that due to the existence of Eps, neither the Loop nor Catmull-Clark subdivision spaces can achieve the optimal approximation rate like the spline spaces on regular regions. Our numerical experiments verify this theory in both Catmull–Clark and Loop subdivisions. In particular, the smaller the valence of the Eps in the mesh, the higher the approximation rate of the corresponding subdivision space projectors. For the Eps of valence $n=3$, the optimal approximation order is achieved in the $L_2$ norm for the Catmull–Clark subdivision, and even in the $L_\infty$-norm for the Loop subdivision.

\begin{figure}[htbp]
	\centering
 \subfigure[Valence $n=5$]{\includegraphics[scale=0.26]{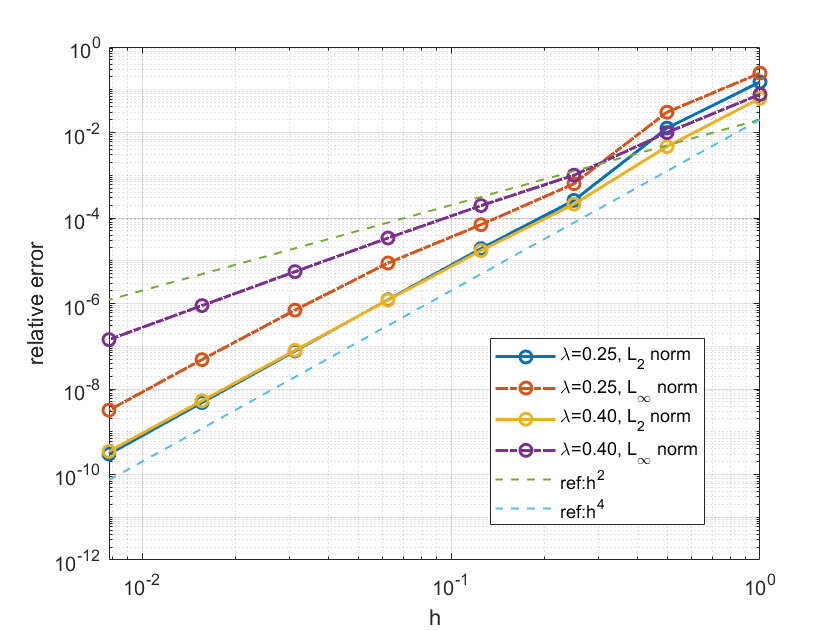}}
\subfigure[Valence $n=8$]{\includegraphics[scale=0.26]{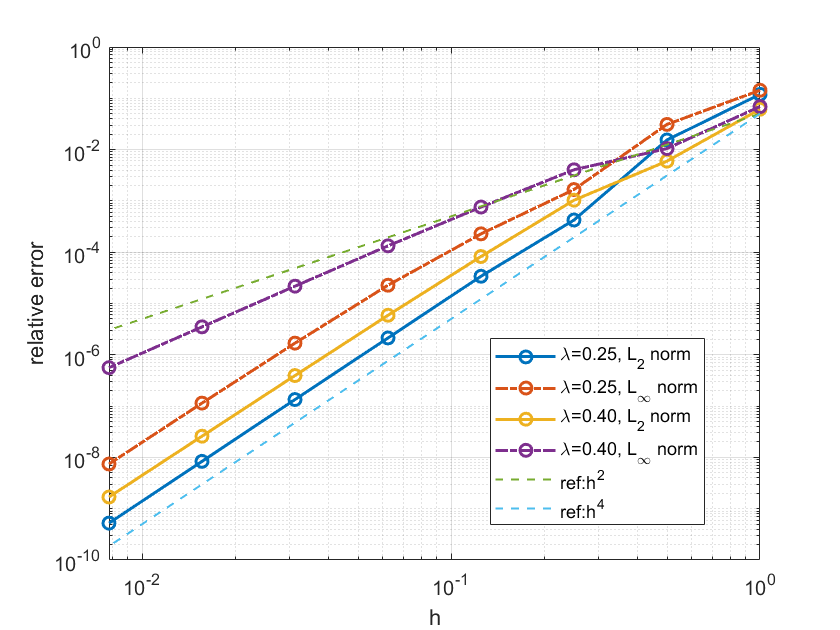}}
	\caption{Relative $L_2$-norm errors and relative $L_{\infty}$-norm errors from the quasi-interpolation projectors for the modified Loop subdivision on the subdivision surfaces with different valances.}
	\label{Fig_20}
\end{figure}

Considering the relationship between the valence $n$ and the subdominant eigenvalue $\lambda$, this inspires us to improve the approximation rate of subdivision surfaces by reducing $\lambda$. The modified Loop subdivision is proposed based on this idea. The subdivision rule of the Eps is controlled by the subdominant eigenvalue. We test the modified Loop subdivision projectors with different subdominant eigenvalues. Fig.\ref{Fig_19} shows the approximation plots of the modified Loop subdivision with Eps of valance $n=5$ and $8$. Numerical experiments show that the optimal approximation order can be achieved under the $L_\infty$-norm when $\lambda=0.25$ and under the $L_2$-norm when $\lambda=0.40$. 


\section{Conclusion \label{conclu}}
In this paper, we propose a general recipe for constructing quasi-interpolation for subdivision surfaces based on point functionals. The existence of extraordinary points makes it challenging to construct a local interpolation problem within a specified domain. Moreover, providing a solution with explicit expressions for the interpolation problem becomes difficult.
In this paper, we constructed the local interpolation problem based on the subdivision matrix and the limit position matrix to tackle the difficulties. We provide the explicit expressions of the quasi-interpolation for the commonly used Loop subdivision and Catmull-Clark subdivision. Numerical experiments indicate that the quasi-interpolation projectors for these two kinds of subdivisions have the approximation order of $3$ under the $L_2$-norm and $2$ under the $L_{\infty}$-norm. In addition, we consider the projector of the modified loop subdivision which has optimal convergence rates in isogeometric analysis. By reducing the subdominant eigenvalues, the modified Loop subdivision projector achieves the optimal approximation order on the mesh with Eps.

Currently, it has not been proved that the proposed quasi-interpolation has the same approximation order as the subdivision space under $L_2$ norm, despite being a projector of the subdivision space. This will be completed in the future. 
The projector characteristic of the proposed quasi-interpolation indicates that this type of approximation can not attain the desired optimal order of approximation (4 for both Loop and Catmull-Clark subdivisions). Moving forward, it is valuable to explore alternative approaches for constructing higher-order approximations. Furthermore, higher-order approximations for subdivision surfaces will be more intriguing.


\bmhead{Acknowledge}
The authors thank the reviewers for providing useful comments and suggestions.

\section*{Declarations}
\bmhead{Conflict of interest} The authors declare no competing interests.

\bibliographystyle{unsrt}
\bibliography{mybibfile}

\begin{thebibliography}{10}

\bibitem{Derose}
T.~Derose, Michael Kass, and Tien Truong.
\newblock Subdivision surfaces in character animation.
\newblock In {\em Proceedings of the 25th Annual Conference on Computer Graphics and Interactive Techniques}, volume~32, pages 85--94, 1998.

\bibitem{dyn_levin_2002}
Nira Dyn and David Levin.
\newblock Subdivision schemes in geometric modelling.
\newblock {\em Acta Numerica}, 11:73--144, 2002.

\bibitem{liao2017subdivision}
Wenhe Liao, Hao Liu, and Tao Li.
\newblock {\em Subdivision Surface Modeling Technology}.
\newblock Springer, 2017.

\bibitem{peters2008subdivision}
J{\"o}rg Peters and Ulrich Reif.
\newblock {\em Subdivision surfaces}.
\newblock Springer, 2008.

\bibitem{DOO1978356}
D.~Doo and M.~Sabin.
\newblock Behaviour of recursive division surfaces near extraordinary points.
\newblock {\em Computer-Aided Design}, 10(6):356--360, 1978.

\bibitem{CATMULL1978350}
E.~Catmull and J.~Clark.
\newblock Recursively generated b-spline surfaces on arbitrary topological meshes.
\newblock {\em Computer-Aided Design}, 10:350--355, 1978.

\bibitem{Loop}
Charles Loop.
\newblock Smooth subdivision surfaces based on triangles.
\newblock Master's thesis, University of Utah, 1987.

\bibitem{Thinshell}
Fehmi Cirak, Michael Ortiz, and Peter Schr{\"o}der.
\newblock Subdivision surfaces: a new paradigm for thin-shell finite-element analysis.
\newblock {\em International Journal for Numerical Methods in Engineering}, 47(12):2039--2072, 2000.

\bibitem{CIRAK2002137}
Fehmi Cirak, Michael~J. Scott, Erik~K. Antonsson, Michael Ortiz, and Peter Schröder.
\newblock Integrated modeling, finite-element analysis, and engineering design for thin-shell structures using subdivision.
\newblock {\em Computer-Aided Design}, 34(2):137--148, 2002.

\bibitem{KANGmodifiedloop}
Hongmei Kang, Wenkai Hu, Zhiguo Yong, and Xin Li.
\newblock Isogeometric analysis based on modified {L}oop subdivision surface with improved convergence rates.
\newblock {\em Computer Methods in Applied Mechanics and Engineering}, 398:115258, 2022.

\bibitem{igasubdthesis}
P.~Barendrecht.
\newblock Isogeometric analysis for subdivision surfaces.
\newblock Master's thesis, Eindhoven University of Technology, 2013.

\bibitem{XIE2020101867}
Jin Xie, Jinlan Xu, Zhenyu Dong, Gang Xu, Chongyang Deng, Bernard Mourrain, and Yongjie~Jessica Zhang.
\newblock Interpolatory {Catmull-Clark} volumetric subdivision over unstructured hexahedral meshes for modeling and simulation applications.
\newblock {\em Computer Aided Geometric Design}, 80:101867, 2020.

\bibitem{Halstead1993}
Mark Halstead, Michael Kass, and Tony DeRose.
\newblock Efficient, fair interpolation using {Catmull-Clark} surfaces.
\newblock In {\em Proceedings of the 20th Annual Conference on Computer Graphics and Interactive Techniques}, pages 35--44, 1993.

\bibitem{WeiyinMa}
Weiyin Ma and Nailiang Zhao.
\newblock {Catmull-Clark} surface fitting for reverse engineering applications.
\newblock In {\em Proceedings of the Geometric Modeling and Processing 2000}, pages 274--283, 2000.

\bibitem{Ma2002}
Weiyin Ma, Xiaohu Ma, Shiu-Kit Tso, and Zhigeng Pan.
\newblock Subdivision surface fitting from a dense triangle mesh.
\newblock In {\em Proceedings of the Geometric Modeling and Processing 2002.}, pages 94--103, 2002.

\bibitem{MA2004525}
Weiyin Ma, Xiaohu Ma, Shiu-Kit Tso, and Zhigeng Pan.
\newblock A direct approach for subdivision surface fitting from a dense triangle mesh.
\newblock {\em Computer-Aided Design}, 36(6):525--536, 2004.

\bibitem{localapprox}
Tom Lyche and Larry~L. Schumaker.
\newblock Local spline approximation methods.
\newblock {\em Journal of {A}pproximation {T}heory}, 15(4):294--325, 1975.

\bibitem{quasimulti}
Carl de~Boor.
\newblock Quasi-interpolants and approximation power of multivariate splines.
\newblock In {\em Computation of Curves and Surfaces}, pages 313--345. Springer, 1990.

\bibitem{deBoor73}
Carl de~Boor and G.~J. Fix.
\newblock Spline approximation by quasi-interpolants.
\newblock {\em Journal of {A}pproximation {T}heory}, 8(1):19--45, 1973.

\bibitem{examples}
Byung-Gook Lee, Tom Lyche, and Knut M{\o}rken.
\newblock Some examples of quasi-interpolants constructed from local spline projectors.
\newblock {\em Mathematical Methods for Curves and Surfaces: {O}slo 2000}, pages 243--252, 2000.

\bibitem{hierarquasi}
Hendrik Speleers.
\newblock Hierarchical spline spaces: quasi-interpolants and local approximation estimates.
\newblock {\em Advances in {C}omputational {M}athematics}, 43(2):235--255, 2017.

\bibitem{effortless}
Hendrik Speleers and Carla Manni.
\newblock Effortless quasi-interpolation in hierarchical spaces.
\newblock {\em Numerische {M}athematik}, 132(1):155--184, 2016.

\bibitem{thbproj}
Alessandro Giust, Bert J{\"u}ttler, and Angelos Mantzaflaris.
\newblock Local {(T)HB}-spline projectors via restricted hierarchical spline fitting.
\newblock {\em Computer {A}ided {G}eometric {D}esign}, 80:101865, 2020.

\bibitem{LRQI}
Francesco Patrizi, Carla Manni, Francesca Pelosi, and Hendrik Speleers.
\newblock Adaptive refinement with locally linearly independent {LR} {B}-splines: Theory and applications.
\newblock {\em Computer Methods in Applied Mechanics and Engineering}, 369:113230, 2020.

\bibitem{KANG2022102147}
Hongmei Kang, Zhiguo Yong, and Xin Li.
\newblock Quasi-interpolation for analysis-suitable {T}-splines.
\newblock {\em Computer Aided Geometric Design}, 98:102147, 2022.

\bibitem{LYCHE2008416}
Tom Lyche, Carla Manni, and Paul Sablonnière.
\newblock Quasi-interpolation projectors for box splines.
\newblock {\em Journal of Computational and Applied Mathematics}, 221(2):416--429, 2008.

\bibitem{2001Approximation}
Greg Arden.
\newblock {\em Approximation properties of subdivision surfaces}.
\newblock PhD thesis, University of Washington, 2001.

\bibitem{2001Fitting}
Nathan Litke, Adi Levin, and Peter Schröder.
\newblock Fitting subdivision surfaces.
\newblock In {\em Proceedings of the Conference on Visualization '01}, pages 319--324, 2001.

\bibitem{dualast}
L.~{Beirao da Veiga}, Annalisa Buffa, D.~Cho, and G.~Sangalli.
\newblock Analysis-suitable {T}-splines are dual-compatible.
\newblock {\em Computer {M}ethods in {A}pplied {M}echanics and {E}ngineering}, 249:42--51, 2012.

\bibitem{subdivisionlinearind}
J\"{o}rg Peters and Xiaobin Wu.
\newblock On the local linear independence of generalized subdivision functions.
\newblock {\em SIAM Journal on Numerical Analysis}, 44(6):2389--2407, 2006.

\bibitem{Boundedcurvature}
Charles Loop.
\newblock Bounded curvature triangle mesh subdivision with the convex hull property.
\newblock {\em The Visual Computer}, 18:316--325, 08 2002.

\end{thebibliography}

\section*{Appendix A.}
\label{appendixB}
\begin{figure}[htbp]
	\centering
 \subfigure[Local Mesh (C-C)]{\includegraphics[scale=0.22]{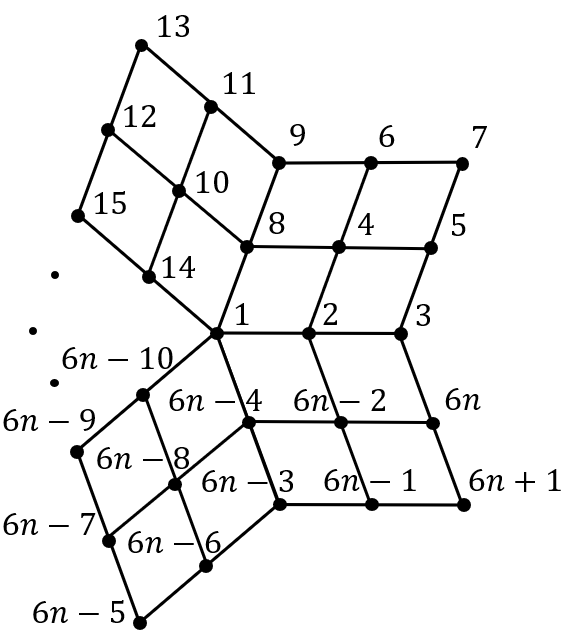}}
 \subfigure[Mesh after subdivision (C-C)]{\includegraphics[scale=0.16]{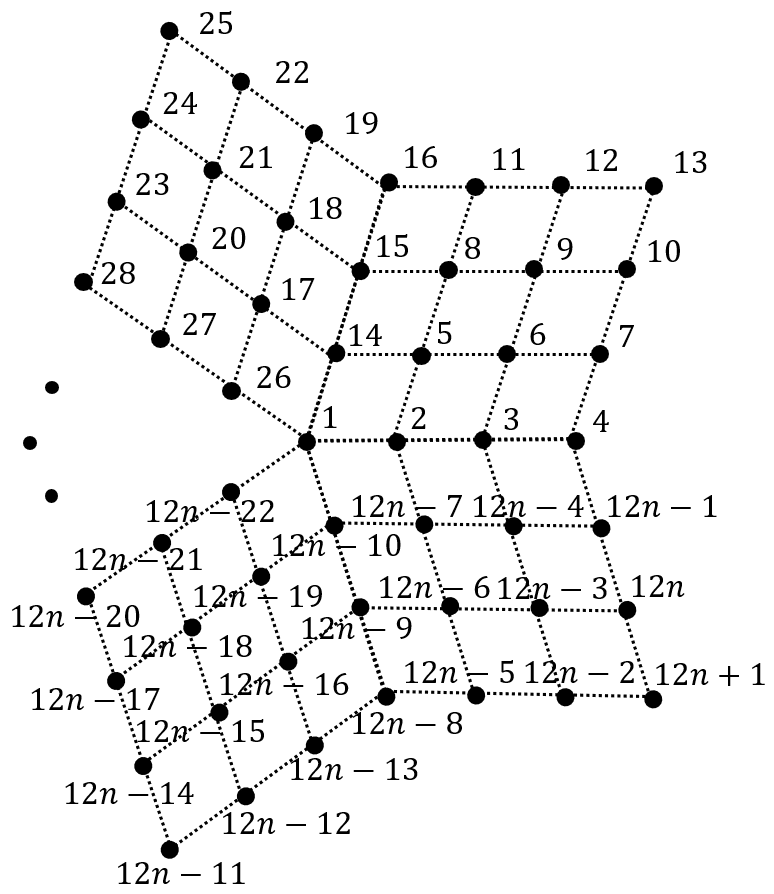}}
 \subfigure[Local Mesh (Loop)]{\includegraphics[scale=0.19]{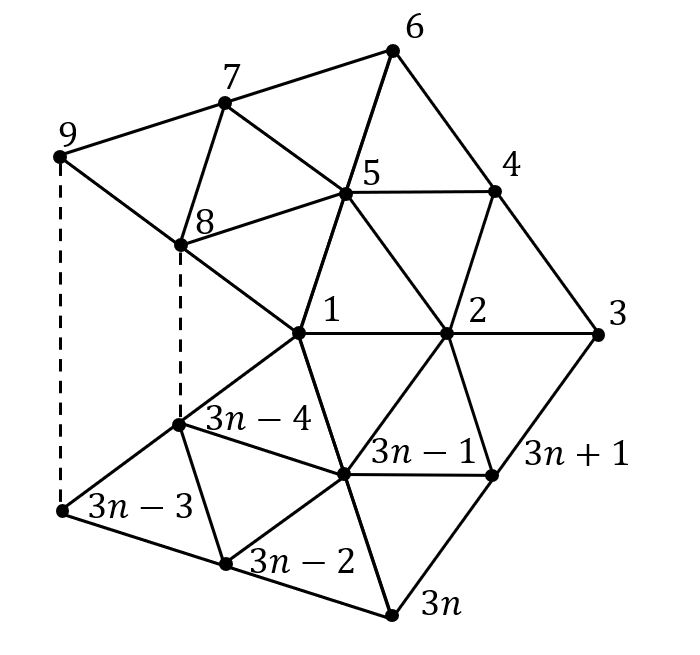}}
 \subfigure[Mesh after subdivision (Loop)]{\includegraphics[scale=0.16]{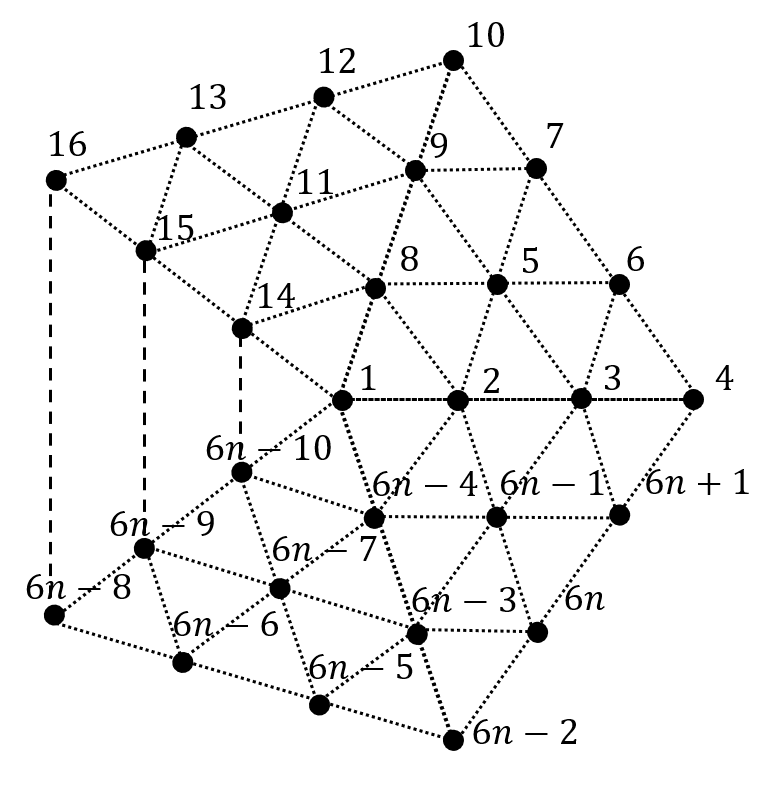}}
	\caption{\label{Fig_3} The numbering scheme for control points, where schemes in (a) and (c) are also applicable to interpolation points.}
\end{figure}

When we number the control points and interpolation points according to Fig.\ref{Fig_3}, the assembled matrices $\bm S$ and $\bm L$ have a clear block structure as follows

\begin{equation}\boldsymbol{S}=\begin{bmatrix}
s_0  &\boldsymbol{s}_1  &\boldsymbol{s}_1  &\dots   &\boldsymbol{s}_1 \\
\boldsymbol{s}_2  &\boldsymbol{S}_0  &\boldsymbol{S}_1  &  &\boldsymbol{S}_2 \\
\boldsymbol{s}_2  &\boldsymbol{S}_2  &\boldsymbol{S}_0  &\ddots  & \\
\vdots   &  &\ddots   &\ddots  &\boldsymbol{S}_1 \\
\boldsymbol{s}_2  &\boldsymbol{S}_1  &  &\boldsymbol{S}_2  &\boldsymbol{S}_0
\end{bmatrix},\quad
\boldsymbol{L}=\begin{bmatrix}
l_0  &\boldsymbol{l}_1  &\boldsymbol{l}_1  &\dots   &\boldsymbol{l}_1 \\
\boldsymbol{l}_2  &\boldsymbol{L}_0  &\boldsymbol{L}_1  &  &\boldsymbol{L}_2 \\
\boldsymbol{l}_2  &\boldsymbol{L}_2  &\boldsymbol{L}_0  &\ddots  & \\
\vdots   &  &\ddots   &\ddots  &\boldsymbol{L}_1 \\
\boldsymbol{l}_2  &\boldsymbol{L}_1  &  &\boldsymbol{L}_2  &\boldsymbol{L}_0
\end{bmatrix} ,   
\end{equation}

\noindent where
\[s_0=1-\frac{7}{4n},\;\boldsymbol{s}_1=\frac{1}{4n^2}\begin{bmatrix}6&0&1&0&0&0\end{bmatrix},\;\boldsymbol{s}_2^T=\frac{1}{64}\begin{bmatrix}24&6&0&16&4&0&4&1&0&0&0&0\end{bmatrix},
\]
\[
\boldsymbol{S}_0=\frac{1}{64}\begin{bmatrix}24&0&4&0&0&0\\36&6&6&1&0&0\\24&24&4&4&0&0\\16&0&16&0&0&0\\24&4&24&4&0&0\\16&16&16&16&0&0\\4&0&24&0&4&0\\6&1&36&6&6&1\\4&4&24&24&4&4\\0&0&16&0&16&0\\0&0&24&4&24&4\\0&0&16&16&16&16\end{bmatrix},\;
\boldsymbol{S}_1=\frac{1}{64}\begin{bmatrix}4&0&0&0&0&0\\1&0&0&0&0&0\\0&0&0&0&0&0\\16&0&0&0&0&0\\4&0&0&0&0&0\\0&0&0&0&0&0\\24&4&0&0&0&0\\6&1&0&0&0&0\\0&0&0&0&0&0\\16&16&0&0&0&0\\4&4&0&0&0&0\\0&0&0&0&0&0\end{bmatrix},\;
\boldsymbol{S}_2=\frac{1}{64}\begin{bmatrix}4&0&4&0&0&0\\1&0&6&0&1&0\\0&0&4&0&4&0\\0&0&0&0&0&0\\0&0&0&0&0&0\\0&0&0&0&0&0\\0&0&0&0&0&0\\0&0&0&0&0&0\\0&0&0&0&0&0\\0&0&0&0&0&0\\0&0&0&0&0&0\\0&0&0&0&0&0\end{bmatrix}.\]
and 
\[l_0=\frac{n}{n+5},\;\boldsymbol{l}_1=\frac{1}{n(n+5)}\begin{bmatrix}4&0&0&1&0&0&0&0&0&0&0&0\end{bmatrix},\;\boldsymbol{l}_2^T=\frac{1}{36}\begin{bmatrix}4&0&1&0&0&0\end{bmatrix},
\]
\[
\boldsymbol{L}_0^T=\frac{1}{36}\begin{bmatrix}
16&4&4&1&0&0\\4&16&1&4&0&0\\0&4&0&1&0&0\\4&1&16&4&4&1\\1&4&4&16&1&4\\0&1&0&4&0&1\\0&0&4&1&16&4\\0&0&1&4&4&16\\0&0&0&1&0&4\\0&0&0&0&4&1\\0&0&0&0&1&4\\0&0&0&0&0&1\end{bmatrix},\;
\boldsymbol{L}_1^T=\frac{1}{36}\begin{bmatrix}1&0&4&0&1&0\\0&0&1&0&4&0\\0&0&0&0&1&0\\0&0&0&0&0&0\\0&0&0&0&0&0\\0&0&0&0&0&0\\0&0&0&0&0&0\\0&0&0&0&0&0\\0&0&0&0&0&0\\0&0&0&0&0&0\\0&0&0&0&0&0\\0&0&0&0&0&0\end{bmatrix}
,\;
\boldsymbol{L}_2^T=\frac{1}{36}\begin{bmatrix}1&0&0&0&0&0\\0&0&0&0&0&0\\0&0&0&0&0&0\\4&1&0&0&0&0\\0&0&0&0&0&0\\0&0&0&0&0&0\\1&4&0&0&0&0\\0&0&0&0&0&0\\0&0&0&0&0&0\\0&1&0&0&0&0\\0&0&0&0&0&0\\0&0&0&0&0&0\end{bmatrix}.
\]
for case of Catmull-Clark subdivision, and 
\[s_0=\alpha,\;\boldsymbol{s}_1=\begin{bmatrix}\beta&0&0\end{bmatrix},\;
\boldsymbol{s}_2^T=\frac{1}{16}\begin{bmatrix}6&1&0&2&0&0\end{bmatrix},\]
\[
\boldsymbol{S}_0=\frac{1}{16}\begin{bmatrix}6&0&0\\10&1&1\\6&6&2\\6&0&2\\6&2&6\\2&0&6\end{bmatrix},\;
\boldsymbol{S}_1=\frac{1}{16}\begin{bmatrix}2&0&0\\1&0&0\\0&0&0\\6&0&0\\2&0&0\\6&2&0\end{bmatrix},\;
\boldsymbol{S}_2=\frac{1}{16}\begin{bmatrix}2&0&0\\1&0&1\\0&0&2\\0&0&0\\0&0&0\\0&0&0\end{bmatrix}.
\]
and
\[l_0=1-n\tau,\;\boldsymbol{l}_1=\begin{bmatrix}\tau&0&0&0&0&0\end{bmatrix},\;
\boldsymbol{l}_2^T=\frac{1}{12}\begin{bmatrix}1&0&0\end{bmatrix},\]
\[
\boldsymbol{L}_0^T=\frac{1}{12}\begin{bmatrix}6&1&1\\1&6&1\\0&1&0\\1&1&6\\0&1&1\\0&0&1\end{bmatrix},\;
\boldsymbol{L}_1^T=\frac{1}{12}\begin{bmatrix}1&0&1\\0&0&1\\0&0&0\\0&0&0\\0&0&0\\0&0&0\end{bmatrix},\;
\boldsymbol{L}_2^T=\frac{1}{12}\begin{bmatrix}1&0&0\\0&0&0\\0&0&0\\1&1&0\\0&0&0\\0&1&0\end{bmatrix}.
\]
for case of Loop subdivision. Therefore, we obtain $\bm{A=LS}$ for any valence $n$.

For case of modified Loop subdivision, the obtained matrices $\bm S$ and $\bm L$ are different from the matrices of Loop subdivision only in the 1-ring neighborhood of the extraordinary point. The corresponding numbers are $1, 2, 5,\cdots,3n-1$ in Fig.\ref{Fig_3}(c) and $1, 2, 8,\cdots,6n-4$ in Fig.\ref{Fig_3}(d). According to the analysis in Section \ref{s.loopGQI}, we obtain different parts $\bm S_n$ and $\bm L_n$ as follow

\[\bm S_{n}=\left[\begin{array}{c|c}
\alpha & \beta \quad \cdots \quad \beta\\
\hline \gamma & \\
\vdots & \hat{\bm S}_{n} \\
\gamma &
\end{array}\right],\quad \bm L_{n}=\left[\begin{array}{c|c}
1-n\tau & \tau \quad \cdots \quad \tau\\
\hline \frac{1}{12}\left(\gamma+\frac58\right) & \\
\vdots & \hat{\bm L}_{n} \\
\frac{1}{12}\left(\gamma+\frac58\right) &
\end{array}\right],\]
where
$\widehat{\bm L}_{n}=circ\left(\frac{1}{12}\left[\gamma_{0}+\frac{45}{8},\gamma_{1}+\frac{7}{8},\gamma_{2},\cdots,\gamma_{n-2},\gamma_{n-1}+\frac{7}{8}\right]\right)$ and 
$\widehat{\bm S}_{n}=circ\left(\left[\gamma_{0},\gamma_{1},\gamma_{2},\cdots,\gamma_{n-2},\gamma_{n-1}\right]\right).$

\begin{figure}[htbp]
	\centering
	\subfigure[Uniform bicubic B-splines]{\includegraphics[scale=0.4]{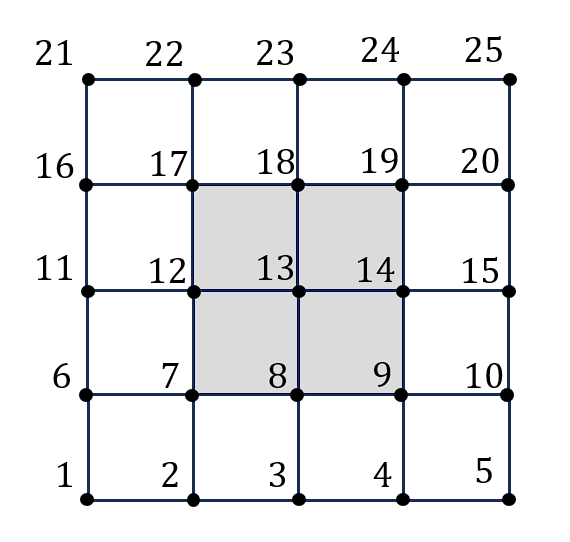}}
	\subfigure[$C^2$-quartic box splines]{\includegraphics[scale=0.4]{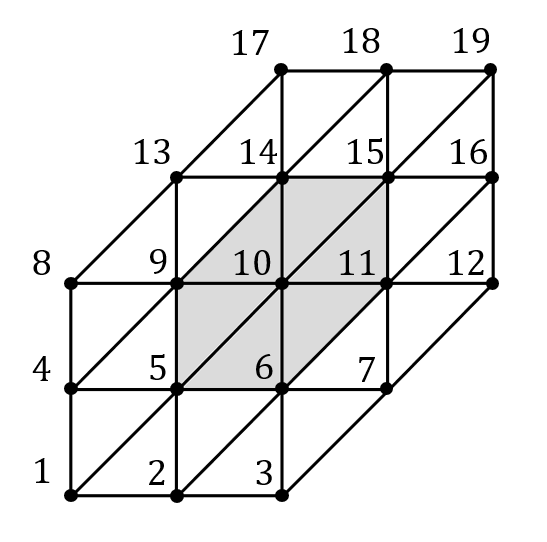}}
	\caption{\label{Fig_12} The numbering scheme for control points, which is also applicable to interpolation points.}
\end{figure}

\section*{Appendix B.}
Our goal is to solve the linear functional $\lambda_i(f)$ for regular points, which is used for treatment in Section \ref{Sec4}.
In this case, the limit surfaces derived from Catmull–Clark subdivision and Loop subdivision are equivalent to bicubic B-spline surfaces and $C^2$ quartic box spline surfaces, respectively. Before solving, we numbered the control points and interpolation points, as shown in Fig.\ref{Fig_12}, where the gray area represents the selected local region. 

Benefiting from the fact that the basis function of bicubic B-splines is a tensor product of the basis function of cubic B-splines, we can easily obtain matrix $\bm A$ through tensor product, i.e
\[\bm A=\frac{1}{48}\begin{bmatrix}
	8 &32  &8  &  & \\
	1&23  &23  &1  & \\
	&8  &32  &8  & \\
	&1  &23  &23  &1 \\
	&  &8  &32 &8
\end{bmatrix}\otimes\frac{1}{48}\begin{bmatrix}
	8 &32  &8  &  & \\
	1&23  &23  &1  & \\
	&8  &32  &8  & \\
	&1  &23  &23  &1 \\
	&  &8  &32 &8
\end{bmatrix},\]
where the operator $\otimes$ represents the Kronecker product, and the resulting $\bm A$ is a matrix of $25\times25$. 

\begin{figure}[htbp]
	\centering
	\subfigure[Operator for point numbered 1]{\includegraphics[scale=0.26]{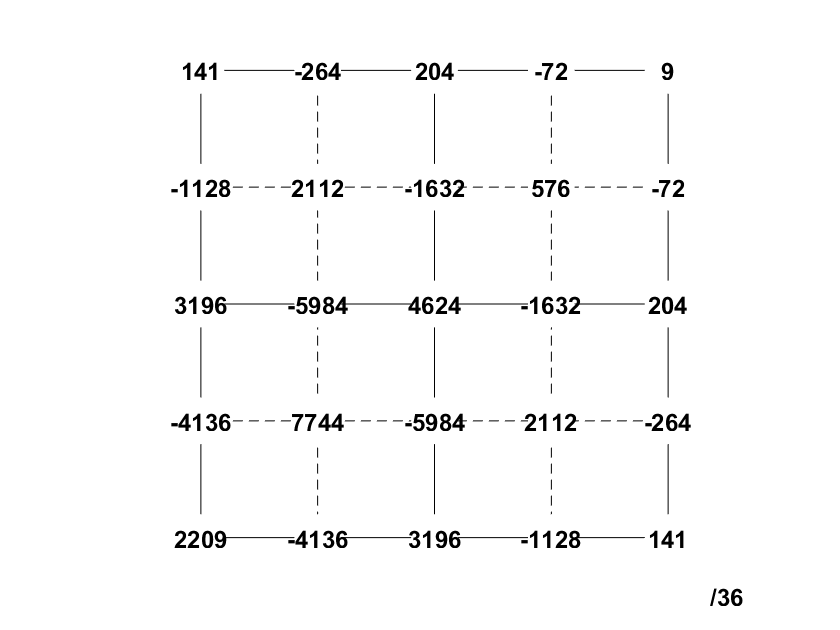}}
	\subfigure[Operator for point numbered 2]{\includegraphics[scale=0.26]{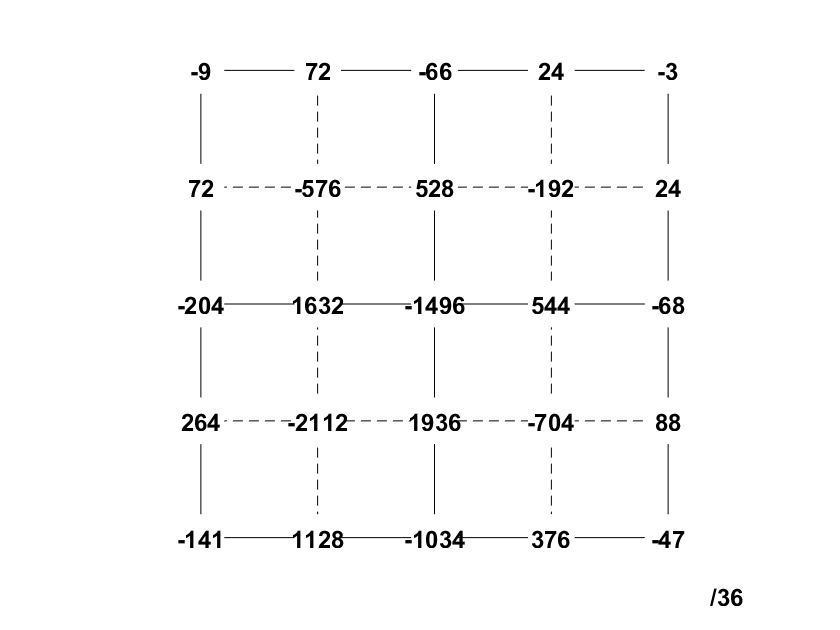}}
	\subfigure[Operator for point numbered 3]{\includegraphics[scale=0.26]{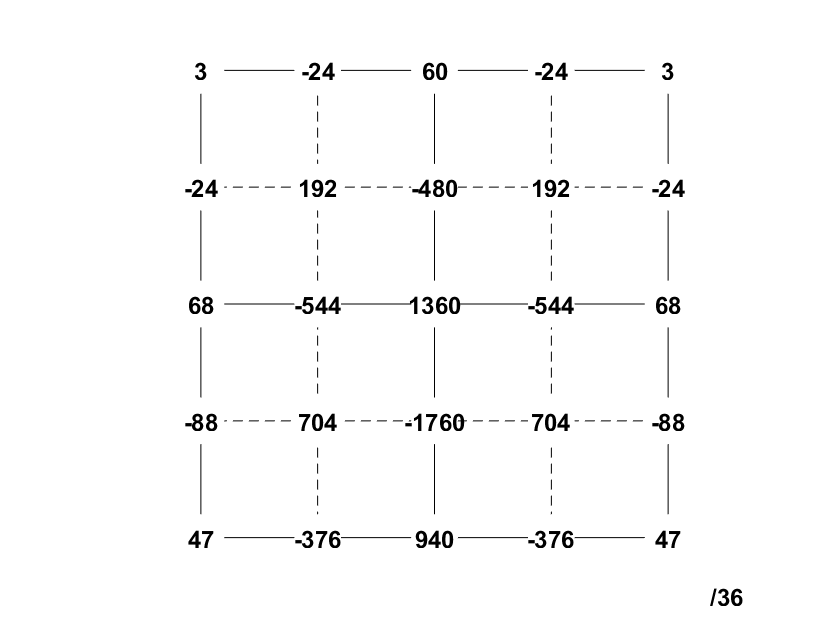}}
	\subfigure[Operator for point numbered 6]{\includegraphics[scale=0.26]{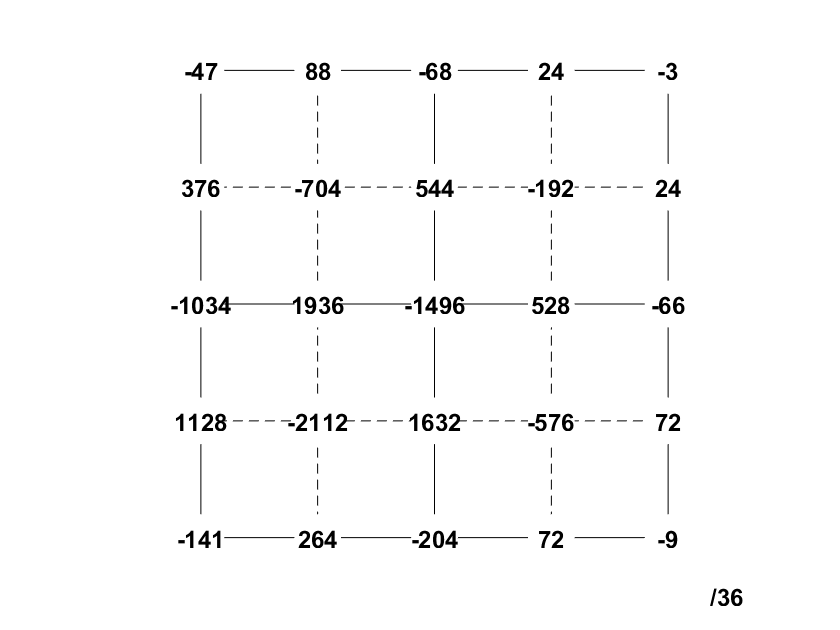}}
 \subfigure[Operator for point numbered 7]{\includegraphics[scale=0.26]{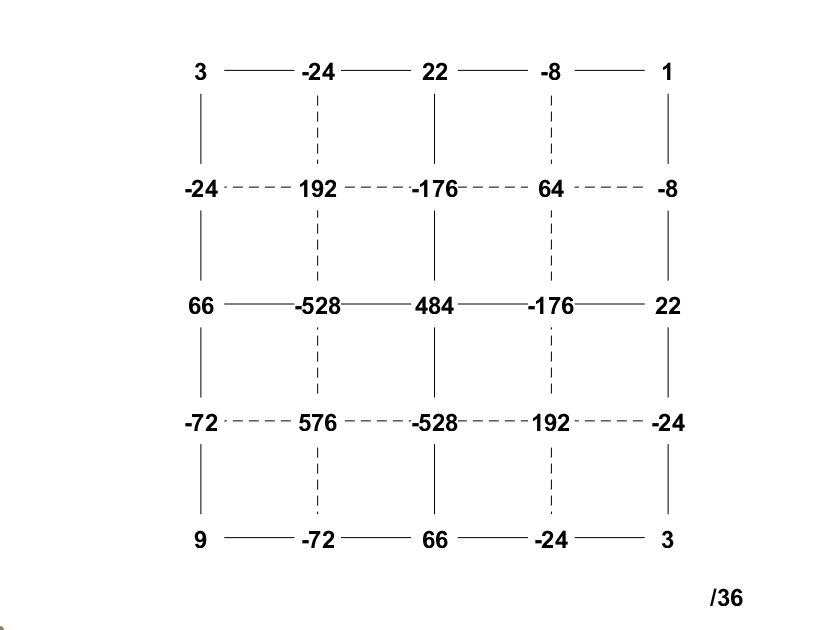}}
 \subfigure[Operator for point numbered 8]{\includegraphics[scale=0.26]{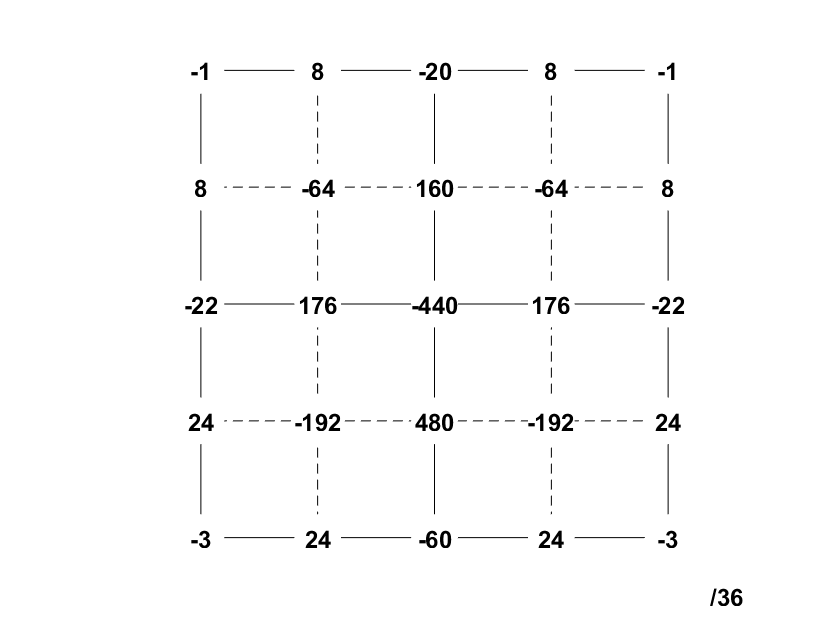}}
	\caption{\label{Fig_11} Linear functionals for uniform bicubic B-splines.}
\end{figure}

The $k$-th row of the inverse of matrix $\bm A$ is coefficient vector $\bm\omega$ of linear functionals for the control point numbered $k$. The operator for the point numbered 13 has been shown as Fig.\ref{F:quasicoeffccc}, while Fig.\ref{Fig_11} shows the quasi interpolation operator for the control points numbered 1, 2, 3, 6, 7, and 8. The operator for the remaining points can be obtained by rotation.

\begin{figure}[htbp]
	\centering
	\subfigure[vertex points]{\includegraphics[scale=0.38]{loopposition6.png}}
	\subfigure[edge points]{\includegraphics[scale=0.19]{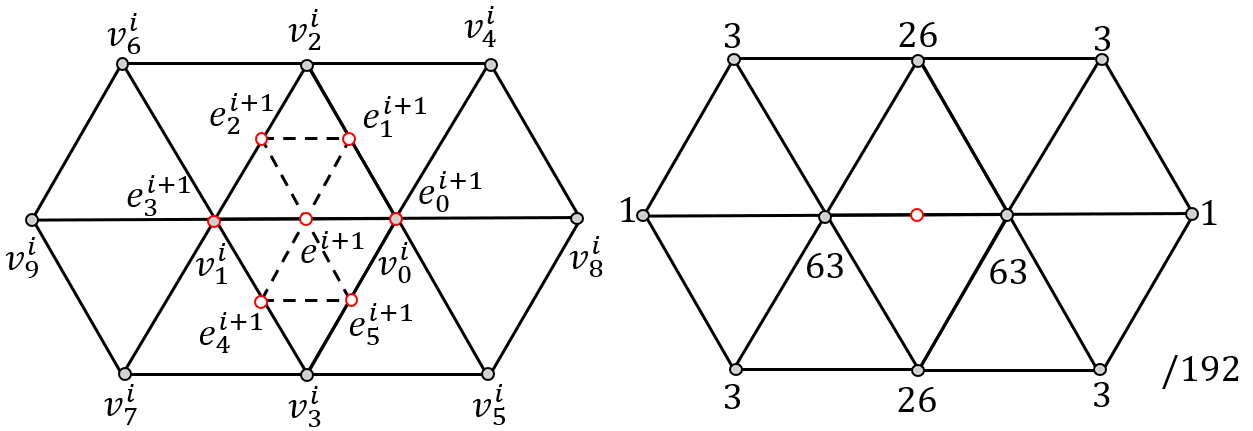}}
	\caption{Position mask for vertex points and edge points.}
	\label{Fig_2}
\end{figure}

For $C^2$ quartic box splines, we can multiply the position matrix by the subdivision matrix to obtain matrix $\bm A$.
This is essentially calculating the formula for the limit position at the interpolation point. There are only two types of interpolation points: vertex and edge points, and limit position formulas are 
\[v^{\infty}=\frac{1}{2}v^i+\frac{1}{12}\sum_{j=0}^{5}e^i_{j},\quad e^{\infty}=\frac{1}{2}e^{i+1}+\frac{1}{12}\sum_{j=0}^{5}e^{i+1}_{j}=\sum_{j=0}^9\varphi_jv_j^i,\]
where $\bm \varphi=\frac{1}{192}[63\;63\;26\;26\;3\;3\;3\;3\;1\;1]$, as shown in Fig.\ref{Fig_2}.

\[
\boldsymbol{A}=\frac{1}{192}\begin{bmatrix}
16&16&0&16&96&16&0&0&16&16&0&0&0&0&0&0&0&0&0\\
3&26&3&1&63&63&1&0&3&26&3&0&0&0&0&0&0&0&0\\
0&16&16&0&16&96&16&0&0&16&16&0&0&0&0&0&0&0&0\\
3&1&0&26&63&3&0&3&63&26&0&0&1&3&0&0&0&0&0\\
1&3&0&3&63&26&0&0&26&63&3&0&0&3&1&0&0&0&0\\
0&3&1&0&26&63&3&0&3&63&26&0&0&1&3&0&0&0&0\\
0&1&3&0&3&63&26&0&0&26&63&3&0&0&3&1&0&0&0\\
0&0&0&16&16&0&0&16&96&16&0&0&16&16&0&0&0&0&0\\
0&0&0&3&26&3&0&1&63&63&1&0&3&26&3&0&0&0&0\\
0&0&0&0&16&16&0&0&16&96&16&0&0&16&16&0&0&0&0\\
0&0&0&0&3&26&3&0&1&63&63&1&0&3&26&3&0&0&0\\
0&0&0&0&0&16&16&0&0&16&96&16&0&0&16&16&0&0&0\\
0&0&0&1&3&0&0&3&63&26&0&0&26&63&3&0&3&1&0\\
0&0&0&0&3&1&0&0&26&63&3&0&3&63&26&0&1&3&0\\
0&0&0&0&1&3&0&0&3&63&26&0&0&26&63&3&0&3&1\\
0&0&0&0&0&3&1&0&0&26&63&3&0&3&63&26&0&1&3\\
0&0&0&0&0&0&0&0&16&16&0&0&16&96&16&0&16&16&0\\
0&0&0&0&0&0&0&0&3&26&3&0&1&63&63&1&3&26&3\\
0&0&0&0&0&0&0&0&0&16&16&0&0&16&96&16&0&16&16
\end{bmatrix}
\]

\begin{figure}[htbp]
	\centering
	\subfigure[Operator for point numbered 1]{\includegraphics[scale=0.26]{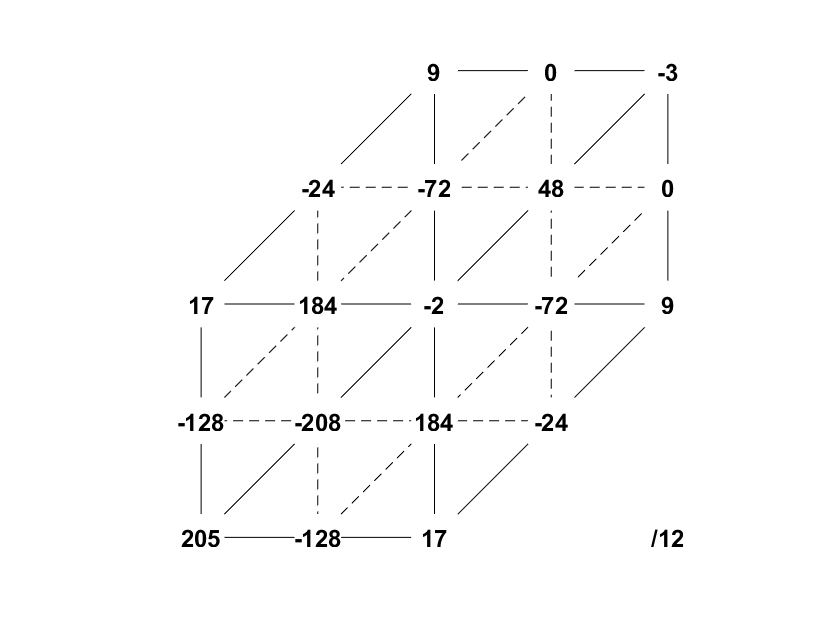}}
	\subfigure[Operator for point numbered 2]{\includegraphics[scale=0.26]{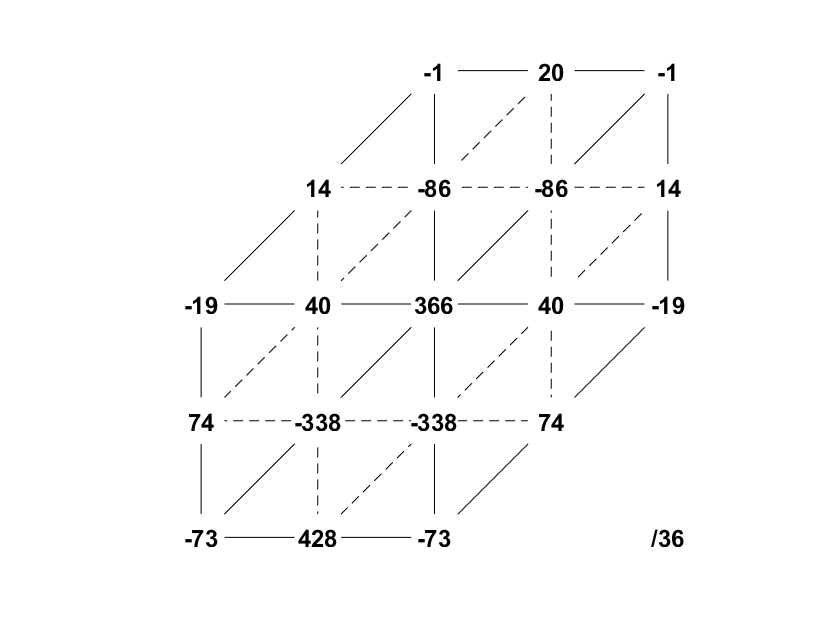}}
	\subfigure[Operator for point numbered 5]{\includegraphics[scale=0.26]{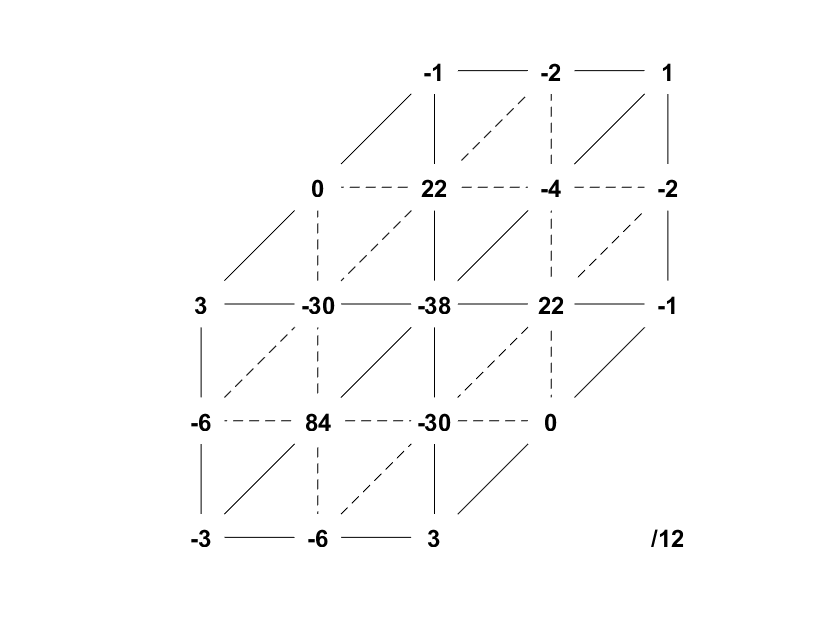}}
 	\subfigure[Operator for point numbered 13]{\includegraphics[scale=0.26]{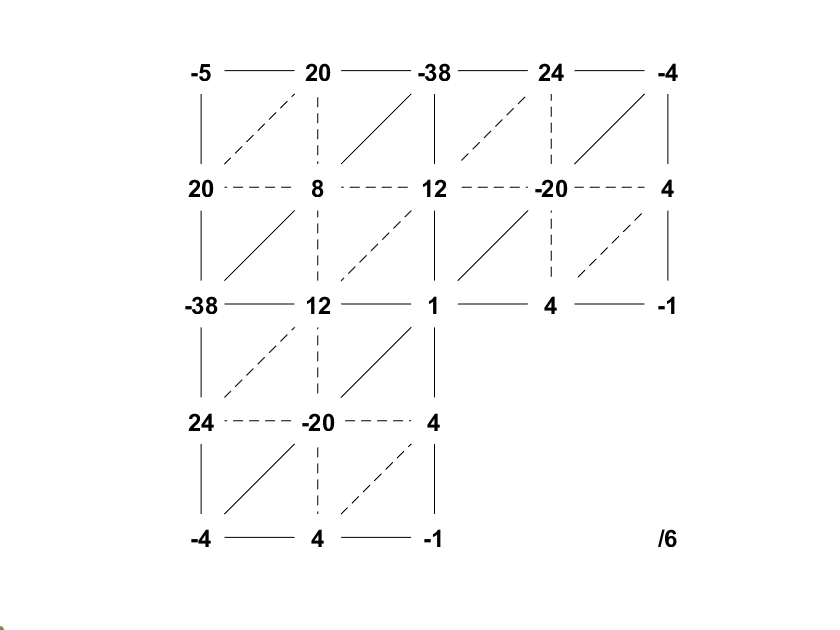}}
	\subfigure[Operator for point numbered 20]{\includegraphics[scale=0.26]{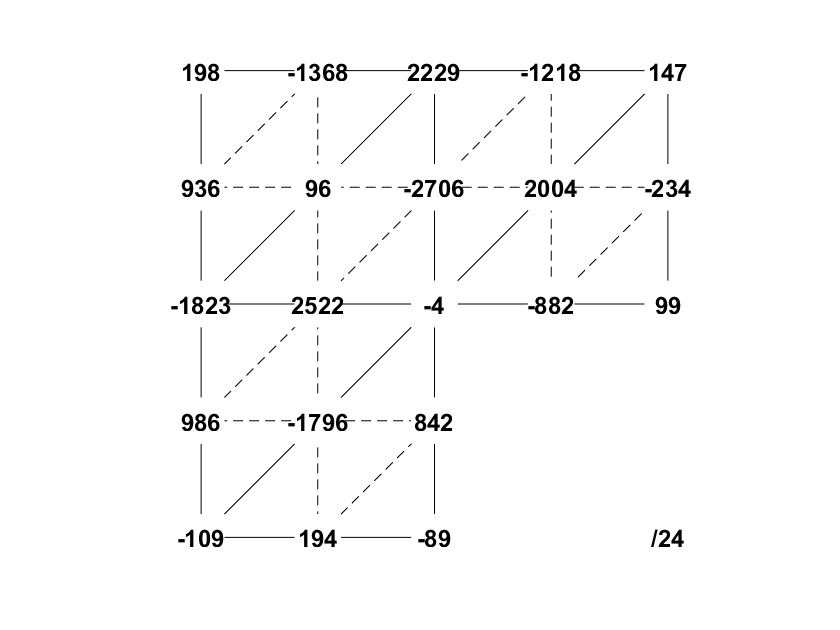}}
	\subfigure[Operator for point numbered 21]{\includegraphics[scale=0.26]{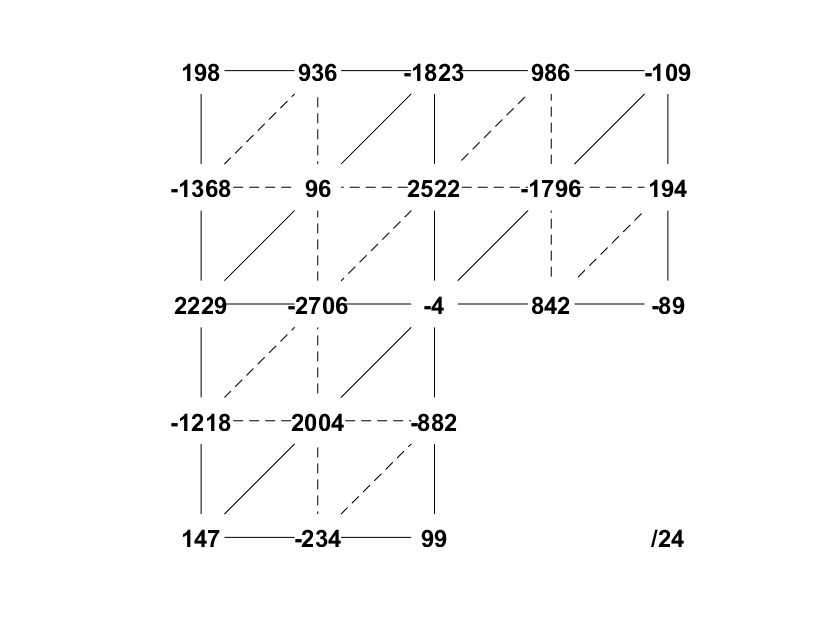}}
	\caption{\label{Fig_13} Linear functionals for $C^2$-quartic box splines. Note that a, b, c and d, e, f are solutions to two different local interpolation problems.}
\end{figure}

Alternatively, we can construct matrix $\bm A$ by directly calculating the value of the basis function at the interpolation point, as \cite{LYCHE2008416} did.
Similarly, invert $\bm A$ to obtain linear functional for each control point. The operator for the point numbered 10 has been shown as Fig.\ref{F:quasicoeffloop}, while Fig.\ref{Fig_13} shows the quasi interpolation operator for the control points numbered 1, 2, and 5. The operator for the remaining points can be obtained by rotation.

\begin{figure}[htbp]
	\centering
	\subfigure[Number of control points]{\includegraphics[scale=0.2]{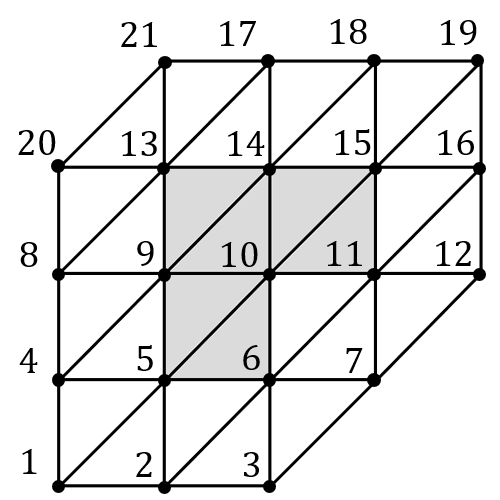}}
	\subfigure[Number of interpolation points]{\includegraphics[scale=0.2]{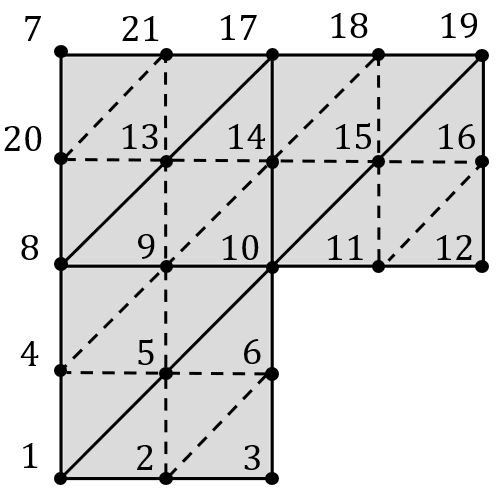}}
	\caption{\label{Fig_14} The numbering scheme for control points and interpolation points in the upper left corner of the triangular mesh.}
\end{figure}

Note that the difference between the mesh in Fig.\ref{Fig_12}(b) and Fig.\ref{Fig_10}(b) is only a few points. Modify Fig.\ref{Fig_12}(b) slightly to obtain the numbering scheme for the control points and interpolation points in the upper left corner of the triangular mesh, as shown in Fig\ref{Fig_14}. 
According to the position masks in Fig.\ref{Fig_2}, we can partially modify the above $\bm A$ to obtain a new matrix as follows for solving the quasi interpolation in the upper left corner of the triangular mesh. 

\[
\boldsymbol{A}=\frac{1}{192}\begin{bmatrix}
16&16&0&16&96&16&0&0&16&16&0&0&0&0&0&0&0&0&0&0&0\\
3&26&3&1&63&63&1&0&3&26&3&0&0&0&0&0&0&0&0&0&0\\
0&16&16&0&16&96&16&0&0&16&16&0&0&0&0&0&0&0&0&0&0\\
3&1&0&26&63&3&0&3&63&26&0&0&1&3&0&0&0&0&0&0&0\\
1&3&0&3&63&26&0&0&26&63&3&0&0&3&1&0&0&0&0&0&0\\
0&3&1&0&26&63&3&0&3&63&26&0&0&1&3&0&0&0&0&0&0\\
0&0&0&0&0&0&0&16&16&0&0&0&96&16&0&0&16&0&0&16&16\\
0&0&0&16&16&0&0&16&96&16&0&0&16&16&0&0&0&0&0&0&0\\
0&0&0&3&26&3&0&1&63&63&1&0&3&26&3&0&0&0&0&0&0\\
0&0&0&0&16&16&0&0&16&96&16&0&0&16&16&0&0&0&0&0&0\\
0&0&0&0&3&26&3&0&1&63&63&1&0&3&26&3&0&0&0&0&0\\
0&0&0&0&0&16&16&0&0&16&96&16&0&0&16&16&0&0&0&0&0\\
0&0&0&1&3&0&0&3&63&26&0&0&26&63&3&0&3&1&0&0&0\\
0&0&0&0&3&1&0&0&26&63&3&0&3&63&26&0&1&3&0&0&0\\
0&0&0&0&1&3&0&0&3&63&26&0&0&26&63&3&0&3&1&0&0\\
0&0&0&0&0&3&1&0&0&26&63&3&0&3&63&26&0&1&3&0&0\\
0&0&0&0&0&0&0&0&16&16&0&0&16&96&16&0&16&16&0&0&0\\
0&0&0&0&0&0&0&0&3&26&3&0&1&63&63&1&3&26&3&0&0\\
0&0&0&0&0&0&0&0&0&16&16&0&0&16&96&16&0&16&16&0&0\\
0&0&0&3&1&0&0&26&63&3&0&0&63&26&0&0&3&0&0&3&1\\
0&0&0&0&0&0&0&3&26&3&0&0&63&63&1&0&26&3&0&1&3\\
\end{bmatrix}
\]

Fig.\ref{Fig_13} shows the linear functionals for the control points numbered 13, 20, and 21. Here we list the operator for control point numberd 13 because the new operator has a smaller $L^1$ norm compared to the original one.

\end{document}